\documentclass[reqno]{amsart}
\usepackage{hyperref}
\usepackage{float}
\usepackage{amsmath}
\usepackage{graphicx} 
\usepackage{latexsym}
\usepackage{amsfonts}
\usepackage{amssymb}
\usepackage{bbold}
\usepackage{amsmath}
\usepackage{xcolor}
\usepackage{dsfont}
\theoremstyle{plain}
\newtheorem{theorem}{Theorem}

\newtheorem{lemma}{Lemma}[section]
\newtheorem{proposition}[lemma]{Proposition}
\theoremstyle{definition}

\newtheorem{definition}[lemma]{Definition}

\newtheorem{remark}[theorem]{Remark}
\newtheorem*{remark*}{Remark}

\title[Expansions for conditioned random walks]{Asymptotic expansions for random walks conditioned to stay positive with arbitrary endpoints}
\author[Tarasov]{Alexander Tarasov}
\address{Hildesheim University, Germany}
\email{tarasov@uni-hildesheim.de}

\newcommand{\pr}{\mathbb{P}}

\newcommand{\lf}{\lfloor}
\newcommand{\rf}{\rfloor}
\newcommand{\E}{\mathbb{E}}

\newcommand{\Z}{\mathbb{Z}}
\newcommand{\cH}{\mathcal{H}}

\begin{document}

\begin{abstract}
We consider a one-dimensional random walk $S_n$ with i.i.d.
integer-valued increments having zero mean and finite variance.
Under suitable moment assumptions, we derive asymptotic expansions
for the local probabilities $\pr(x+S_n=y,\tau_x>n)$, where
$\tau_x=\inf\{n\ge1:x+S_n\le0\}$.
These extend previous expansions for walks started at zero
to arbitrary starting points. We obtain different expansions
according to the sizes of the starting point and endpoint relative
to $\sqrt n$.
\end{abstract}

\maketitle

\section{Introduction}
Let $\{X_k\}$ be a sequence of independent, identically distributed
integer-valued random variables with maximal span $1$, zero mean
$\E X_1=0$ and finite variance $\sigma^2:=\E X_1^2>0$.
We consider the random walk $\{S_n; n\ge0\}$ defined by $S_0=0$ and
\begin{align*}
    S_n:=X_1+X_2+\ldots+X_n,\ n\ge1.
\end{align*}
Consider the stopping time $\tau_x$ defined as
$$
\tau_x:=\inf\{n\ge1:x+S_n\le0\},\quad x\in \Z_+ := \{ 0,1,2 \dots\}.
$$
By a random walk conditioned to stay positive we mean the walk $\{S_k\}_{k=0}^n$ 
on the event $\{\tau_x > n\}$. It is known that 
$\pr(\tfrac{S_n}{\sigma\sqrt{n}} \in \cdot\mid \tau_0>n)$ 
converges to a Rayleigh distribution, see \cite{Caravenna05}. On a functional level the walk 
converges to a Brownian meander, see \cite{Iglehart74,Bolthausen76}.

The goal of the present paper is to obtain asymptotic expansions 
for a local probability for random walks conditioned to stay positive 
with arbitrary endpoints:
$$
    \pr(x+S_n = y, \tau_x > n), \quad x,y \in \Z_+.
$$
This extends the expansions for walks started at zero obtained in 
our joint work with Denisov and Wachtel~\cite{DTW23,DTW25a}.

For walks started at zero, local limit theorems were proved
by Caravenna~\cite{Caravenna05} in the domain of attraction
of the normal law and by Vatutin and Wachtel~\cite{WV09}
in the $\alpha$-stable setting.
For positive starting points, Uchiyama~\cite{Uchiyama2011}
obtained local asymptotics for lattice walks with finite variance,
while Doney~\cite{Doney2012} treated walks in the domain
of attraction of a stable law.

An interesting connection of such walks to other areas of probability 
theory can be found in Green and Sawhney's 
work~\cite{GreenSawhney2026} on the limiting probability that a random permutation has an 
invariant set of size $k$. 
They note that the local estimate they need follows
from~\cite{DTW25a} and their appendices give a 
shorter proof for their particular walk, adopting  
the Berry--Esseen arguments that were developed in~\cite{DTW24,DTW_BE2}.
Barrier options in financial mathematics are another possible application. 
Their pricing in discrete-time models leads to random walks conditioned 
to stay positive \cite{BGK97,LinPalmer2013}.
Such walks also arise in counting paths on $\Z^2$ 
that go only up and right while staying above a line; see \cite{WB}. 
The question is open for lines with irrational slope. 
This model leads to non-integer-valued walks, which we plan to consider in future work.
The current paper can be seen as an important step towards analysis of general walks.

\smallskip

For random walks conditioned to stay positive, different regimes 
appear naturally. It is already visible in the survival 
probability $\pr(\tau_x>n)$. For fixed $x$, it decays as 
$n^{-1/2}$, with a coefficient depending on $x$ and the 
increment distribution. When $x \sim u \sqrt n, u > 0$,
the survival probability instead converges to a positive
limit given by Brownian motion (see, for example, \cite[Section~1]{Doney2012}).
For local probabilities, the scales of both endpoints matter, see Doney~\cite{Doney2012}. 
At leading order, Grama and Xiao~\cite{GramaXiaoLocal25} obtained a 
heat-kernel approximation valid across 
all three regimes, with an explicit error estimate.

\subsection{Navigation of the paper}
In the next section we formulate the main results, compare them with the existing literature and 
discuss our method.
Section~\ref{sec:2} contains the proof of
Theorem~\ref{thm:1}. Section~\ref{sec:exitestimates}
gives preliminary estimates used in the proofs of
Theorems~\ref{thm:2} and~\ref{thm:3}, presented in
Sections~\ref{sec:3} and~\ref{sec:4}, respectively.

\section{Main results}

The form of the expansion in the following theorems depends on the 
scale of the endpoints
relative to $\sqrt n$. The remainder bounds below hold uniformly over the 
stated ranges, but the form of the remainders determines where each expansion 
gives a better approximation.
We say that the point $z$ is near the boundary if it is $o(\sqrt{n})$. In such
regimes the variable $z$ appears in the numerator of the remainder. 
The first theorem deals with the regime
where both $x$ and $y$ are near the boundary.
\begin{theorem}
\label{thm:1}
Assume that $\E|X_1|^{r+3}$ is finite for some integer $r\ge 0$. 
Assume also that $S_n$ is integer-valued with maximal span $1$.
Then there exist functions
$U_1$, $U_2$, \ldots, $U_{\lf \frac{r+2}{2}\rf}$ such that, 
uniformly in $y \ge 1$, $x\ge 1$,
\begin{multline} \label{eq:thm1}
    \pr(x+S_n=y, \tau_x>n)
=
    \sum_{j=1}^{\lf \frac{r+2}{2}\rf}
        \frac{U_j(x,y)}{n^{j+1/2}}\\
+
    O\left(\frac{(1+\min(x,y)) (x+y)^{r+2}}{n^{(r+4)/2}}
    \log^{{\lf \frac{r+2}{2}\rf}}n\right).
\end{multline}
Moreover,
\begin{align*}
    |U_j(x,y)|
\le
    C(1+\min(x,y)) (x+y)^{2j-1}.
\end{align*}
\end{theorem}

To describe the leading coefficient in \eqref{eq:thm1}, put
\[
    V(x):=-\E S_{\tau_x}=x+\E[-x-S_{\tau_x}],\qquad x\ge0,
\]
and define $V^*$ in the same way for the walk $S^*_n =-X_1 - \dots -X_n$. 
These functions are harmonic
for the walks killed upon leaving $(0,\infty)$:
\[
    V(x)=\E[V(x+X_1);\,x+X_1>0],\qquad x>0.
\]
With this normalization, the classical first-passage asymptotics read
\begin{equation}\label{eq:survivalasymptotics}
    \pr(\tau_x>n)\sim
    \sqrt{\frac{2}{\pi}}\frac{V(x)}{\sigma\sqrt n},
    \qquad n\to\infty,
\end{equation}
for every fixed $x\ge0$; see, for example, \cite{DSW18}.
Using \eqref{eq:udefinition} from the proof of
Theorem~\ref{thm:1}, we obtain
\[
    U_1(x,y)=\frac{2V(x)V^*(y)}{\sigma^3\sqrt{2\pi}}.
\]
This coincides with the results of
Doney~\cite{Doney2012}.

Decomposing the probability $\pr(x+S_{n+1} = y, \tau_x > n+1)$ according to the position
after the first step or before the last step gives recursions in $x$
and $y$, respectively. Comparing these recursions term by term with
\eqref{eq:thm1} shows that, on the positive half-line, $U_j$ is
polyharmonic of order at most $j$ in each variable for the killed walk
and its dual. Polyharmonic coefficients in expansions for small-step walks in the
quadrant were studied by Chapon, Fusy and
Raschel~\cite[Section~3]{ChaponFusyRaschel2020} and
Nessmann~\cite[Theorems~4.1 and~4.2]{Nessmann2024}.

For a walk started at zero, Theorem 1.5 in \cite{DTW23} gives under the assumptions of Theorem~\ref{thm:1}
\begin{align*}
    \pr(S_n = y, \tau_0 > n)
=
    \sum_{j=1}^{\big\lf \frac{r+2}{2} \big \rf}
        \frac{W_j(y)}{n^{j+1/2}}
+
    O \left( \frac{(1+y)^{r+2}}{n^{(r+4)/2}} \log^{\big\lf \frac{r+1}{2} \big \rf} n \right)
\end{align*}
for some polyharmonic functions $W_j$ with $W_1 = cV^*$ for some constant $c$. So Theorem~\ref{thm:1} 
is a natural generalization of Theorem 1.5 in \cite{DTW23}.

\smallskip

The next expansion extends Theorem 1 in~\cite{DTW25a}, which treats
walks started at zero, to positive starting points. It is designed
for starting points near the boundary and endpoints of order
$\sqrt n$. Reversing time gives the corresponding expansion with
the roles of the two endpoints interchanged.
\begin{theorem} \label{thm:2}
Assume that $\E|X_1|^{r+3}$ is finite for some integer $r\ge0$. 
Assume also that $S_n$ is integer-valued with maximal span $1$.
Then there exist functions $P_2(t,x), P_3(t,x), \ldots, P_{r+2}(t,x)$ which are, 
for every fixed $x$, polynomials in $t$, such that, uniformly in $x, y \ge 1$,
\begin{align*}
    \pr(x+S_n=y,\tau_x>n)
&=
    e^{-\frac{y^2}{2\sigma^2n}}
    \sum_{\nu=2}^{r+2}
        \frac{P_\nu\bigl(\tfrac{y}{\sigma\sqrt{n}},x\bigr)}{n^{\nu/2}}\\
&\hspace{1.8cm}+
    o\left(\frac{1}{(\sqrt{n}+y)^{r+2}}\right)
+
    O\left(\frac{x^{r+2}\log^{\lf \frac{r+2}{2} \rf+1} n}{(\min(\sqrt{n},y))^{r+3}}\right).
\end{align*}
Moreover,
\[
    P_2(t,x)=t\,\frac{2V(x)}{\sigma^2\sqrt{2\pi}},
\]
$P_\nu$ has degree at most $3\nu-5$ and contains 
$t^i$ only for $i \equiv \nu+1 \pmod 2$
and, for every $x\ge1$ and $t\in\mathbb R$,
\begin{align*}
    |P_\nu(t,x)| \le C(1+|t|)^{3\nu-5}x^{\nu-1}.
\end{align*}
\end{theorem}
The functions $P_\nu(t,x)$ can be calculated explicitly via moments of $X_1$ and the sums
\begin{align*}
    \sum_{u\ge1}\sum_{z\ge0}
        z^h U_{\ell+1}(x,u)
        \pr(X_1=-u-z),
\end{align*}
where $U_\ell$ are functions from Theorem~\ref{thm:1}.

The leading term agrees with the local asymptotics obtained by
Uchiyama~\cite[Theorem~1.3]{Uchiyama2011} and
Doney~\cite[Proposition~11]{Doney2012}. 

\smallskip

Finally, the next theorem covers the regime where both $x$ and $y$ 
are of order $\sqrt{n}$.
\begin{theorem} \label{thm:3}
Assume that $\E|X_1|^{r+3}$ is finite for some integer $r\ge1$.
Assume also that $S_n$ is integer-valued with maximal span $1$.
Then there exist polynomials
$B_1, \ldots, B_{r+1}$ and $ q_1, \ldots, q_{r+1}$ with
$B_1 \equiv q_1 \equiv 1/(\sigma\sqrt{2\pi})$, such that,
uniformly in $x,y\ge1$,
\begin{align*}
&    \pr(x+S_n=y,\tau_x>n)\\
&\hspace{1cm}=
    e^{-\frac{(y-x)^2}{2\sigma^2n}}
    \sum_{\nu=1}^{r+1}
        \frac{q_\nu\bigl(\tfrac{y-x}{\sigma\sqrt{n}}\bigr)}{n^{\nu/2}}
-
    e^{-\frac{(x+y)^2}{2\sigma^2n}}
    \sum_{\nu=1}^{r+1}
        \frac{B_\nu\bigl(\tfrac{x}{\sigma\sqrt{n}},\,
                           \tfrac{y}{\sigma\sqrt{n}}\bigr)}{n^{\nu/2}}\\
&\hspace{8cm}
+
    o\left(\frac{n}{\min(\sqrt{n},x,y)^{r+3}}\right).
\end{align*}
\end{theorem}
The polynomials $q_\nu$ are the coefficients in the local
expansion for the unconditioned walk
\cite[Theorem~VII.3.16]{Petrov}. The polynomials $B_\nu$
are determined by the moments of $X_1$ and the sums
\[
    \sum_{u\ge1}\sum_{z\ge0}
        z^h
        P_\mu^*(t,u)\pr(X_1=-u-z),
\]
where $P_\mu^*$ are the coefficients of Theorem~\ref{thm:2}
for the walk $-S_n$.

Since $B_1\equiv q_1\equiv 1/(\sigma\sqrt{2\pi})$, the leading term
in the expansion is
\[
    \frac{1}{\sigma\sqrt{2\pi n}}
    \left(
        e^{-\frac{(y-x)^2}{2\sigma^2n}}
        -
        e^{-\frac{(y+x)^2}{2\sigma^2n}}
    \right).
\]
This is $1/(\sigma\sqrt n)$ times the
transition density at time $1$ of a standard Brownian motion started
at $x/(\sigma\sqrt n)$ and killed at zero, evaluated at
$y/(\sigma\sqrt n)$; see the reflection principle
in~\cite[Section~2.2]{MP10}.

For lattice walks, Ford~\cite{Ford2009} obtained
$1-e^{-2xy/(\sigma^2 n)}$ as a leading-order approximation to the
probability of staying positive conditional on the endpoint.
Multiplying it by the unconditioned local probability gives
the difference of the two Gaussian kernels appearing above.
The work of Korolyuk~\cite{Korolyuk1962} yields the first three terms
of the local probability considered here when both endpoints are of
order $\sqrt n$.
In the strongly non-lattice case, Siegmund and
Yuh~\cite{SiegmundYuh1982} obtained the first correction to the
Brownian approximation for the corresponding distribution function.

Under exponential moment assumptions, Borovkov~\cite{Borovkov1962} 
obtained asymptotic expansions for the joint distribution of the maximum and 
the value of the walk at time $n$, assuming an absolutely 
continuous component. 
Theorem~\ref{thm:2} can also be used to derive an expansion of $\pr(\tau_x >n)$ 
for $x$ of order $\sqrt{n}$ by first deriving $\pr(\tau_x = n)$ and then summing.
Related expansions for the distribution of the maximum of a random walk 
were obtained by Nagaev~\cite{Nagaev1970} under nonlattice assumptions.

\begin{remark}
The same method should give expansions for
$\pr(x+S_n\ge y,\tau_x>n)$ in the regimes considered above.
When $y$ is near the boundary, use
\[
    \pr(x+S_n\ge y,\tau_x>n)
    =\pr(\tau_x>n)
     -\sum_{z=1}^{y-1}\pr(x+S_n=z,\tau_x>n),
\]
together with appropriate expansions for the survival probability. 
For $x$ near the boundary, such expansions are provided by \cite[Theorem~1.4]{DTW23}.
When $y$ is of order $\sqrt n$, the proofs of
Theorems~\ref{thm:2} and~\ref{thm:3} consider $\pr(x+S_n = y, \tau_x < n)$ and 
split the trajectory at $\tau_x$. With $x+S_{\tau_x} = -z$, the rest of the path 
satisfies $S_{n - \tau_x} =z+ y$. For the probability 
$\pr(x+S_n\ge y,\tau_x>n)$, the corresponding term is $\pr(S_{n - \tau_x} \ge z + y)$ 
instead and one shall use Petrov's expansion for the distribution function
\cite[Theorem~VI.3.6]{Petrov} in place of the local expansion.
\end{remark}

\begin{remark}

Similarly to $\tau_x$ one can define
$$
    \overline{\tau}_x:=\inf\{n\ge1:x+S_n < 0\},\quad x\in \Z_+.
$$
Since the walk is integer-valued, $\overline{\tau}_x=\tau_{x+1}$, hence
\begin{align*}
    \pr(x+S_n = y, \overline{\tau}_x > n) 
&= 
    \pr(x+S_n = y, \tau_{x+1} > n)\\
&=
    \pr((x+1)+S_n = y+1, \tau_{x+1} > n).
\end{align*}
So, expansions for $\overline{\tau}_x$ can be deduced from those for $\tau_x$
and we do not formulate them separately.
\end{remark}

\smallskip

We do not seek optimal moment assumptions or remainder
estimates in Theorems~\ref{thm:2} and~\ref{thm:3}
to avoid additional technical details and believe that one more term
can be added to each expansion, including Theorem 1 in~\cite{DTW25a}.
Our expansions for conditioned walks 
fully rely on Petrov's expansions for unconditional walks; see~\cite{Petrov}.
More accurate decompositions would require 
keeping track of cancellations between the
remainders in Petrov's expansion instead of bounding them in absolute value. 
Similar cancellations based on the reflection principle were used in our
works on the Berry--Esseen inequality~\cite{DTW24,DTW_BE2}.

\subsection{About the method}

The foundation of our approach here is the expansion of local probabilities
obtained for a walk started at zero and ending near the
boundary in~\cite{DTW23}:
\begin{align} \label{eq:methodmain}
    \pr(S_n = y, \tau_0 > n)
=
    \sum_{j=1}^{\big\lf \frac{r+2}{2} \big \rf}
        \frac{W_j(y)}{n^{j+1/2}}
+
    O \left( \frac{(1+y)^{r+2}}{n^{(r+4)/2}} \log^{\big\lf \frac{r+1}{2} \big \rf} n \right)
\end{align}
It is applicable only near the boundary since the remainder 
grows with $y$ and is comparable to the leading term when $y\sim \sqrt{n}$.

\medskip
This expansion happens to be the key to all other regimes.
For both starting and end points near the boundary, one can split the path at its first 
minimum and reverse the time before reaching it
\begin{align*}
    &\pr(x+S_n = y, \tau_x > n)\\
&\hspace{2cm}=
    \sum_{k=0}^n
    \sum_{z = 1}^{\min(x,y)}
        \pr(-S_k = x-z, \tau^*_{0} > k) \pr(S_{n-k} =y-z, \overline{\tau}_0 > n-k),
\end{align*}
where $\tau^*_{0} = \inf\{n\ge 1, -S_n \le 0 \}$. Then both parts of the path 
are walks conditioned 
to stay positive and started at level zero, so one can apply 
the decomposition \eqref{eq:methodmain} and its analog for $\overline{\tau}_0$.
In this way we prove Theorem~\ref{thm:1} and obtain a decomposition with a remainder 
that grows as $(x+y)^{r+3}$.

\medskip

In all other regimes we start with
\begin{align*}
    \pr(x+S_n =y, \tau_x > n) = \pr(S_n = y-x) - \pr(x+S_n = y, \tau_x <n).
\end{align*}
Since expansions of $\pr(S_n = y-x)$ are known, we can consider only the probability $\pr(x+S_n = y, \tau_x <n)$. 
The advantage of this probability is its natural decomposition at time $\tau_x$.
For $x=0$ and $y\approx \sqrt{n}$, we apply the Markov property at time $\tau_0$ and $\tau_0-1$:
\begin{align*}
    &\pr(S_n = y, \tau_0 < n)\\
& \hspace{1.5cm}=
     \sum_{k=1}^n \sum_{z \ge 0}
     \sum_{u\ge 1}
        \pr(S_{k-1} = u, \tau_0 > k-1)
        \pr(X_k = -u-z)
        \pr(-z + S_{n-k} = y)
\end{align*}
The last $n-k$ steps form an unconditional walk, so standard Edgeworth expansions from~\cite{Petrov}
are applicable. To expand $\pr(S_{k-1} = u, \tau_0 > k-1)$ one can apply \eqref{eq:methodmain} with $u$ 
instead of $y$ and 
the remainder, which grows with $u$, is controlled by $\pr(X_k = -u-z)$, since $\E |X_1|^{r+3} < \infty$.
In this way, Theorem 1 in~\cite{DTW25a} is proven, yielding a remainder 
$O\big( \frac{1}{(\min(\sqrt{n}, y))^{r+3}} \big)$. Unlike \eqref{eq:methodmain} it is non-increasing 
in $y$ but provides less accuracy for small $y$.

\medskip

For Theorem~\ref{thm:2}, when $x$ is near the boundary and $y\approx \sqrt{n}$ the 
proof is almost the same. We write
\begin{align} \label{eq:methodmain2}
    &\pr(x+S_n = y, \tau_x < n)\\
& \hspace{1cm}=
     \sum_{k=1}^n \sum_{z \ge 0}
     \sum_{u\ge 1}
        \pr(x+S_{k-1} = u, \tau_x > k-1)
        \pr(X_k = -u-z)
        \pr(-z + S_{n-k} = y).
\end{align}
The only change is that we use Theorem~\ref{thm:1} instead of Theorem 1.5 in~\cite{DTW23} 
to expand $\pr(x+S_{k-1} = u, \tau_x > k-1)$.

Alternatively, in this case one could decompose the path by minimum, obtaining two paths 
started at zero, one ending near the boundary 
and one ending at the scale $\sqrt{n}$, applying Theorem 1.5 from~\cite{DTW23} and Theorem 1 from~\cite{DTW25a}
 respectively. 
It would be interesting to compare the form of the 
coefficients, but that would require another full proof.

\medskip

Finally, for both starting and end points being of order $\sqrt{n}$, we again use \eqref{eq:methodmain2}. 
The difference is that now the walk up to time $\tau_x-1$, 
conditioned to stay positive,
starts on the scale $\sqrt{n}$ and ends near the boundary. In this way one can 
deduce Theorem~\ref{thm:3} from Theorem~\ref{thm:2} with reversed time. In this case splitting by minimum 
does not give the results directly, since the minimum can be large and we have to sum 
over all its possible values. 

\medskip

Another important remark is that all results of the current paper and \cite{DTW23,DTW25a} 
should be valid for random walks with bounded density. Basically, the method converts 
asymptotic decompositions for unconditional walks into decompositions for conditional walks.
Since similar decompositions exist not only for integer-valued walks but also for walks with 
bounded density (see \cite{Petrov}), the method is applicable. Similarly, it is 
possible to generalize the results to the non-integer-valued lattice case, when $X_1$ 
takes values on a shifted lattice
$\delta + \Z$, $\delta \notin \Z$. However, it requires separate treatment and a
modification of the theory developed in \cite{DTW23}, since it involves an oscillation 
in contrast to the case of bounded density or integer-valued walks with maximal span $1$.

\section{Proof of Theorem \ref{thm:1}} \label{sec:2}
Recall from \cite{DTW23} the definition of the class $\mathcal{H}_r^d$.
\begin{definition}
\label{def.1}
Fix some $d\in\mathbb{N}_0$ and $r\in\mathbb{N}_0 \cup \{-1\}$.
Let $\mathcal{H}_r^d$ denote the class of generating functions \begin{equation*}
H(s) = \sum\limits_{n=0}^\infty h_n s^n, \ s \in [0,1)
\end{equation*}
of sequences $\{h_n\}_{n=0}^\infty$ satisfying the following two conditions:
\begin{enumerate}
    \item $h_n = O\left( \frac{(\log n)^d}{n^{(r+3)/2}}\right)$;
    \item $\sum\limits_{n=0}^\infty h_n P(n) = 0$ for every polynomial $P$ such that $\deg P \le \lfloor r/2 \rfloor$.
\end{enumerate}
We will write that
\begin{align*}
    H(s) = G(s) + \cH_r^d
\end{align*}
if the difference belongs to $\cH_r^d$:
\begin{align*}
    H(s) - G(s) \in \cH_r^d.
\end{align*}
\end{definition}

Without loss of generality, we can assume that $x \ge y$. 
The case $x < y$ can be handled in exactly the same way. 
Theorem~\ref{thm:1} follows from the following two propositions.
\begin{proposition} \label{prop:1thm1}
    There are functions $u_j(x,y)$ such that for every fixed pair $x,y$
\begin{align}\label{eq:prop5}
    \sum_{n=0}^\infty 
        s^n \pr(x+S_n = y, \tau_x > n) 
=
    \sum_{i=0}^{r+1} 
        u_i(x,y)(1-s)^{i/2} 
+
    \mathcal{H}_{r+1}^{\lf \frac{r+1}{2} \rf},
\end{align}
Moreover,
\begin{align*}
    |u_i(x,y)| \le Cyx^{i}.
\end{align*}
\end{proposition}
The dependence of the remainder on $x$ and $y$ is given by
the following proposition.
\begin{proposition}\label{prop:2thm1}
There is a constant $C$ such that, for $n \ge 2$,
\begin{align*}
    \left|
        \pr(x+S_n = y, \tau_x > n) 
    -
        \sum_{i=0}^{r+1} 
            u_i(x,y)[s^n](1-s)^{i/2}
    \right|
\le 
    C
    \frac{yx^{r+2}\log^{\lf \frac{r+2}{2}\rf}n}{n^{(r+4)/2}}.
\end{align*}
\end{proposition}

To deduce Theorem~\ref{thm:1}, note that for even $i$
the sequence $[s^n](1-s)^{i/2}$ has finite support.
For the odd powers, the binomial coefficient expansion in
\cite[Section 2]{DTW23} gives
\[
    [s^n](1-s)^{\ell-1/2}
    =\sum_{j=\ell}^{\lf(r+2)/2\rf}
        \frac{\beta_{j+1}^{(\ell+1)}}{n^{j+1/2}}
    +O\bigl(n^{-\lf(r+2)/2\rf-3/2}\bigr),
    \qquad 1\le\ell\le\lf(r+2)/2\rf,
\]
with constants $\beta_{j+1}^{(\ell+1)}$. Thus the coefficients
in \eqref{eq:thm1} are
\[
    U_j(x,y)=\sum_{\ell=1}^j
        \beta_{j+1}^{(\ell+1)}u_{2\ell-1}(x,y).
\]
The bounds in Proposition~\ref{prop:1thm1} give
$|U_j(x,y)|\le Cyx^{2j-1}$ for $x\ge y\ge1$.
The error from this change of basis is
$O(yx^{r+1}n^{-\lf(r+2)/2\rf-3/2})$, which is absorbed
by the remainder in Proposition~\ref{prop:2thm1}.
This proves the theorem.

\begin{proof}[Proof of Proposition \ref{prop:1thm1}]
Let $S^*$ be an independent copy of the reflected walk $-S$,
and write
\[
    \tau^*:=\inf\{n\ge1:S_n^*\le0\}.
\]
Recall also that
$\overline{\tau}:=\overline{\tau}_0=\inf\{n\ge1:S_n<0\}$.
Assuming $\pr(S_0 = 0) = \pr(S^*_0 = 0) = 1$ and splitting the path by its first minimum, we can write
\begin{multline} \label{eq:splitbymin}
    \pr(x+S_n = y, \tau_x >n ) 
\\
=
    \sum_{z = 1}^{x \wedge y}
    \sum_{k = 0}^{n}
        \pr(S^*_k = x-z, \tau^* > k) \pr(S_{n-k} = y-z, \overline{\tau} > n-k).
\end{multline}
Since we restrict ourselves to the case $x\ge y$, we have $x\wedge y = y$.
Multiplying this equality by $s^n$ and summing over $n\ge0$, we obtain
\begin{align} \label{eq:split_by_min_gen}
    \sum_{n=0}^\infty &s^n \pr(x + S_n = y, \tau_x > n)
\nonumber \\
&=
    \sum_{n=0}^\infty 
    \sum_{z=1}^y
    \sum_{k=0}^{n}
        s^n \pr(S_k^* = x - z, \tau^* > k)
        \pr(S_{n-k} = y - z, \overline{\tau} > n-k)\nonumber \\
&=
    \sum_{z=1}^y
    \left(\sum_{n=0}^\infty
        s^n \pr(S_n^* = x - z, \tau^* > n) \right)
    \left(\sum_{n=0}^\infty
        s^n \pr(S_n = y - z, \overline{\tau} > n)\right).
\end{align}
Write $a_n^{(j)}=[s^n](1-s)^{j-3/2}$.
Under the assumptions of Theorem~\ref{thm:1},
\cite[Theorem 1.5]{DTW23} gives
\begin{align}\label{eq:locprobdecomp}
    \pr(S_n=x,\overline{\tau}>n)=\sum_{j=1}^{\lf \frac{r+2}{2}\rf} \overline{U}_j(x)a^{(j+1)}_n
+
    O\left(\frac{(1+x)^{r+2}}{n^{(r+4)/2}}\log^{{\lf \frac{r+1}{2}\rf}}n\right),
\end{align}
for some functions $\overline{U}_j$ that depend on the random walk.
By \cite[Lemma 4.3]{DTW23} we then have
\begin{align} \label{eq:Hdecomplocprob}
    \sum_{n=0}^\infty
        s^n\pr(S_n = x, \overline{\tau} > n)
=
    \sum_{i=0}^{r+1}
        q_i(x)(1-s)^{i/2}
+
    \mathcal{H}_{r+1}^{\lf \frac{r+1}{2}\rf},
\end{align}
where $q_{2j-1} = \overline{U}_j$ and the coefficients $q_{2j}$ are 
determined by the remainder in \eqref{eq:locprobdecomp}, 
see \cite[Lemma 4.3]{DTW23}. Their values are not needed below, 
we use only the bound $|q_i(x)|\le C(1+x)^i$, which holds for all 
$i \le r+1$ by \cite[Lemma 6.2]{DTW23}. That lemma is stated for 
$i\le r$, but its proof gives the same bound for $i=r+1$ without any changes.

Again by \cite[Theorem 1.5]{DTW23} one has
\begin{align} \label{eq:locprobdecompviagenfun*}
    \pr(S^*_n=x,\tau^*>n)=\sum_{j=1}^{\lf \frac{r+2}{2}\rf} U^*_j(x)a^{(j+1)}_n
+
    O\left(\frac{(1+x)^{r+2}}{n^{(r+4)/2}}\log^{{\lf \frac{r+1}{2}\rf}}n\right),
\end{align}
for some functions $U_j^*$. By \cite[Lemma 4.3]{DTW23} we then have
\begin{align*}
    \sum_{n=0}^\infty
        s^n\pr(S^*_n = x, \tau^* > n)
=
    \sum_{i=0}^{r+1}
        q^*_i(x)(1-s)^{i/2}
+
    \mathcal{H}_{r+1}^{\lf \frac{r+1}{2}\rf},
\end{align*}
where $q^*_{2j-1} = U_j^*$.

Recall \cite[Lemma 4.12]{DTW23}: if $H(s) \in \mathcal{H}_{r}^d$ 
for some $r \ge 0$, then for every integer $k \ge 1$
\begin{align*}
    H(s)(1-s)^{k/2} 
\in 
    \mathcal{H}_{r}^d.
\end{align*}
Obviously, it holds also for $k=0$. Using this lemma and \cite[Lemma 4.7]{DTW23},
which states that the class $\mathcal{H}_{r+1}^{\lf \frac{r+1}{2}\rf}$ is closed under multiplication, one can easily get that
\begin{multline} \label{eq:genfunprod1}
    \left(\sum_{n=0}^\infty
        s^n \pr(S_n^* = x, \tau^* > n) \right)
    \left(\sum_{n=0}^\infty
        s^n \pr(S_n = y, \overline{\tau} > n)\right)\\
=
    \left(
        \sum_{i=0}^{r+1}
        q^*_i(x)(1-s)^{i/2}
    \right)
    \left(
        \sum_{j=0}^{r+1}
        q_j(y)(1-s)^{j/2}
    \right)
+
    \mathcal{H}_{r+1}^{\lf \frac{r+1}{2}\rf}.
\end{multline}
By \cite[Lemma 4.10]{DTW23}, for $i+j > r+1$, we have
\begin{align*}
    (1-s)^{(i+j)/2} \in \mathcal{H}_{i+j-1}^0 \subset \mathcal{H}_{r+1}^0.
\end{align*}
Hence we can expand the brackets on the RHS of \eqref{eq:genfunprod1} and obtain
\begin{multline*}
    \left(\sum_{n=0}^\infty
        s^n \pr(S_n^* = x, \tau^* > n) \right)
    \left(\sum_{n=0}^\infty
        s^n \pr(S_n = y, \overline{\tau} > n)\right)\\
=
    \sum_{\ell =0}^{r+1}
        \left(
            \sum_{i = 0}^\ell
                q^*_i(x)q_{\ell-i}(y)
        \right)
        (1-s)^{\ell/2}
+
    \mathcal{H}_{r+1}^{\lf \frac{r+1}{2}\rf}.
\end{multline*}
Substituting this into \eqref{eq:split_by_min_gen},
we obtain \eqref{eq:prop5} for
\begin{align} \label{eq:udefinition}
    u_\ell(x,y)
=
    \sum_{j=0}^{\ell}
        \sum_{z=1}^y
        q_{j}(y-z)q^*_{\ell-j}(x-z).
\end{align}
Finally, by \cite[Lemma 6.2]{DTW23}, we have
\begin{align} \label{eq:qpower}
    |q_i(x)| + |q^*_i(x)| \le C(1+x)^i,
\end{align}
and hence
\begin{align*}
    |u_i(x,y)| \le Cyx^i.
\end{align*}
\end{proof}
\begin{proof} [Proof of Proposition \ref{prop:2thm1}]
By the local limit theorem
(see \cite[Chapter~VII]{Petrov}),
\begin{align} \label{eq:exitLLT}
    \sup_{v\in\mathbb Z}\pr(S_n=v)\le \frac C{\sqrt n},
    \qquad n\ge1.
\end{align}
Write $\tau_y^*=\inf\{k\ge1:y+S_k^*\le0\}$.
For $n\ge2$, let $\ell=\lf n/2\rf$.
Reversing time on the last $n-\ell$ steps and using
\eqref{eq:exitLLT}, we obtain
\begin{align} \label{eq:oneendpointbound}
    \pr(x+S_n=y,\tau_x>n)
&\le \sum_{v\ge1}\pr(x+S_\ell=v)
       \pr(y+S^*_{n-\ell}=v,\tau_y^*>n-\ell)\nonumber\\
&\le \frac C{\sqrt n}\pr(\tau_y^*>n-\ell)
\le \frac{Cy}{n}.
\end{align}
Here $\pr(\tau_y^*>m)\le Cy/\sqrt m$ follows from
the proof of \cite[eq.~(46)]{Doney2012} for $y\le\sqrt m$,
and from the bound by one otherwise.
The case $n=1$ is trivial.

Corollary~12 in \cite{Doney2012} gives the bound below
for $x,y\le\sqrt n$.
For $x\ge\sqrt n$ or $y\ge\sqrt n$, it follows from
\eqref{eq:oneendpointbound}. The case $x,y\ge\sqrt n$ covered by \eqref{eq:exitLLT}.
Thus, for all $x,y\ge1$ and $n\ge1$,
\begin{align} \label{eq:basice_xy_estim}
    \pr(x+S_n = y, \tau_x > n)
\le
    C\frac{xy}{n^{3/2}}.
\end{align}
Using this, we can prove the proposition for $x \ge \sqrt{n}$. In this case we can write
\begin{multline*}
    \big|
        \pr(x+ S_n = y , \tau_x >n)
    -
        \sum_{i=0}^{{r+1}}
            u_i(x,y) [s^n](1-s)^{i/2}
    \big|\\
\le
    \pr(x+ S_n = y , \tau_x >n)
    +
    \sum_{i=0}^{{r+1}}
            |u_i(x,y)[s^n](1-s)^{i/2}|.
\end{multline*}
Then, for $x\ge \sqrt{n}$, \eqref{eq:basice_xy_estim} gives
\begin{align*}
    \pr(x+S_n = y, \tau_x > n)
\le
    C\frac{xy}{n^{3/2}}
\le
    C\frac{yx^{r+2}}{n^{\frac{r+4}{2}}}.
\end{align*}
Since by Proposition~\ref{prop:1thm1} $|u_i(x,y)| \le Cyx^{i}$, we have for odd $i \le {r+1}$
\begin{align*}
        |u_i(x,y) [s^n](1-s)^{i/2}|
\le
    C\frac{yx^{i}}{n^{i/2+1}}
\le
    C\frac{yx^{r+2}}{n^{\frac{r+4}{2}}}.
\end{align*}
For even $i$, the estimate is immediate, since $[s^n](1-s)^{i/2}=0$ for $n>i/2$.

{For the case $x \le \sqrt{n}$, we need more precise estimates. 
We will use the following simple bound. Letting $s\uparrow1$ in 
\eqref{eq:Hdecomplocprob}, we get 
$\sum_{k=0}^\infty \pr(S_k = u, \overline{\tau} > k) = q_0(u)$, 
and $q_0$ is bounded by \cite[Lemma 6.2]{DTW23}. Hence
\begin{align} \label{eq:greenbound1}
    \sum_{z=1}^y
    \sum_{k=0}^{n}
        \pr(S_{k} = y-z, \overline{\tau} > k) \le C(1+y),
\end{align}
and, in the same way,
\begin{align} \label{eq:greenbound2}
    \sum_{z=1}^y
    \sum_{k=0}^{n}
        \pr(S_k^* = x-z, \tau^* > k) \le C(1+y).
\end{align}}
Since we have already proved Proposition~\ref{prop:2thm1} for $x \ge \sqrt{n}$, the following lemma completes its proof and hence the proof of Theorem~\ref{thm:1}.
\end{proof}
\begin{lemma} \label{lem:xlesqrtn}
    For $n\ge 2$ and $x \le \sqrt{n}$, there is a constant $C > 0$ such that
\begin{align*}
    \big|
        \pr(x+ S_n = y , \tau_x >n)
    -
        \sum_{i=0}^{r+1}
            u_i(x,y) [s^n](1-s)^{i/2}
    \big|
\le
    C
    {
        \frac{yx^{{r+2}}}{n^{{\frac{r+4}{2}}}}
    }
    \log^{\lf \frac{r+2}{2}\rf} n.
\end{align*}
\end{lemma}
\begin{proof}
Using \eqref{eq:udefinition} and \eqref{eq:splitbymin}, we obtain
\begin{multline} \label{eq:prob_u_diff}
        \pr(x + S_n = y, \tau_x > n)
    -
        \sum_{i=0}^{r+1}u_i(x,y)[s^n](1-s)^{i/2}\\
=
    \sum_{z = 1}^y
    \sum_{k = 0}^{n}
        \pr(S^*_k = x-z, \tau^* > k) \pr(S_{n-k} = y-z, \overline{\tau} > n-k)\\
    -
        \sum_{z=1}^y
        \sum_{i=1}^{r+1}
        \bigg(
        \sum_{j=0}^{i}
            q_{j}(y-z)q^*_{i-j}(x-z)
        \bigg)
            [s^n](1-s)^{i/2}.
\end{multline}
For $i>r+1$ and $x \le \sqrt{n}$, we have
\begin{align*}
    |q_j(y-z)q{^*}_{i-j}(x-z)[s^n](1-s)^{i/2}|
\le
    C\frac{x^{{i}}}{n^{(i+2)/2}}
\le
    C\frac{x^{r+2}}{n^{(r+4)/2}}.
\end{align*}
Hence
\begin{multline} \label{eq:diff_prodqq*}
    \Bigg|
    \sum_{k=0}^n
    \left(
        \sum_{i=0}^{r+1} q_i({y-z})[s^k](1-s)^{i/2}
    \right)
\cdot
    \left(
        \sum_{i=0}^{r+1} q_i^*({x-z})[s^{n-k}](1-s)^{i/2}
    \right)\\
-
    \sum_{i=0}^{r+1}
    \left(
        \sum_{j=0}^i
            q_j(y-z)q^*_{i-j}(x-z)
    \right)
        [s^n](1-s)^{i/2}
    \Bigg|
\le
    C\frac{x^{r+2}}{n^{(r+4)/2}}.
\end{multline}
For brevity, let
\begin{align*}
    p_n(x) = \pr(S_n=x, \overline{\tau} > n)
,
\quad
\quad
    p_n^*(x) = \pr(S^*_n=x, \tau^* > n),
\end{align*}
and
\begin{align*}
    p_n^{\text{reg}}(x) 
=
    \sum_{i=0}^{r+1}
        q_i(x)[s^n](1-s)^{i/2},
\quad \quad
    p_n^{\text{reg}*}(x) 
=
    \sum_{i=0}^{r+1}
        q^*_i(x)[s^n](1-s)^{i/2}.
\end{align*}
Summing \eqref{eq:diff_prodqq*} over $z=1,\ldots,y$, we obtain
\begin{multline*}
    \Bigg|
        \sum_{z=1}^y
        \Bigg[
        \sum_{k=0}^n
            p_k^{\text{reg}}(y-z)
            p_{n-k}^{\text{reg}{*}}(x-z)\\
-
        \sum_{i=0}^{r+1}
        \left(
            \sum_{j=0}^i
                q_j(y-z)q^*_{i-j}(x-z)
        \right)
            [s^n](1-s)^{i/2}
        \Bigg]
    \Bigg|
\le
    C\frac{yx^{r+2}}{n^{(r+4)/2}}.
\end{multline*}
Combining this with \eqref{eq:prob_u_diff}, we obtain
\begin{multline} \label{eq:lem8:3}
    \left|
        \pr(x + S_n = y, \tau_x > n)
    -
        \sum_{i=0}^{r+1}u_i(x,y)[s^n](1-s)^{i/2}
    \right|
\le 
    C\frac{yx^{r+2}}{n^{(r+4)/2}}\\
+
    \left|
        \sum_{z=1}^{y}
        \bigg[
        \sum_{k=0}^n
            p_k(y-z)p_{n-k}^*(x-z)
    -
        \sum_{k=0}^n
            p_k^{\text{reg}}(y-z)
            p_{n-k}^{\text{reg}{*}}(x-z)
        \bigg]
    \right|.
\end{multline}
Thus, it is enough to prove
\begin{multline} \label{eq:probconvregdiff}
    \left|
        \sum_{z=1}^{y}
        \bigg[
        \sum_{k=0}^n
            p_k(y-z)p_{n-k}^*(x-z)
    -
        \sum_{k=0}^n
            p_k^{\text{reg}}(y-z)
            p_{n-k}^{\text{reg}{*}}(x-z)
        \bigg]
    \right|\\
\le
    C\frac{yx^{r+2}}{n^{(r+4)/2}}\log^{\lf \frac{r+2}{2}\rf} n.
\end{multline}

We bound only the sum for $k$ from $0$ to $\lf n/2 \rf-1$, the sum for $k \ge \lf n/2 \rf$ can be estimated similarly.
For every $u,v$, with $N = \lf n/2 \rf-1$, we have
\begin{align} \label{eq:splitregrest}
    \sum_{k=0}^{N}\nonumber
        p_k(u)
        p_{n-k}^*(v)
-
    \sum_{k=0}^{N}
        p_k^{\text{reg}}(u)
        p_{n-k}^{\text{reg}*}(v)
&=
    \sum_{k=0}^{N}
        p_k(u)
        \big[
            p_{n-k}^*(v)
        -
            p_{n-k}^{\text{reg}*}(v)
        \big]\\
&+
    \sum_{k=0}^{N}
        \big[
            p_{k}(u)
        -
            p_{k}^{\text{reg}}(u)
        \big]
        p_{n-k}^{\text{reg}*}(v).
\end{align}
{By \eqref{eq:greenbound1}}
\begin{align*}
    \sum_{z=1}^y
    \sum_{k=0}^N
        p_{k}(y-z)
\le
    Cy.
\end{align*}
Since $q^*_{2j-1}(x) = U^*_j(x)$, it follows from the decomposition \eqref{eq:locprobdecompviagenfun*} that
\begin{align*}
    |p_{n-k}^*(v)
-
    p_{n-k}^{\text{reg}*}(v)|
\le
    C\frac{x^{r+2}}{n^{(r+4)/2}}\log^{{\lf \frac{r+1}{2}\rf}}n, \quad \quad k\le N.
\end{align*}
Combining the last two estimates, we obtain
\begin{align} \label{eq:lem8firstsum}
    \left|
    \sum_{z=1}^y
    \sum_{k={0}}^{{N}}
        p_{k}(y-z)
        \big[p_{n-k}^{*}(x-z)
    -
        p_{n-k}^{\text{reg}*}(x-z)\big]
    \right|
\le
    C\frac{yx^{r+2}}{n^{(r+4)/2}}\log^{{\lf \frac{r+1}{2}\rf}}n.
\end{align}
By the definition of $p_{n-k}^{\text{reg}}$ and the decomposition \eqref{eq:Hdecomplocprob}, we have
\begin{align*}
    \sum_{n=0}^\infty
    \big(
        p_{n}(x)
    -
        p_{n}^{\text{reg}}(x)
    \big) s^n
\in
    \mathcal{H}_{r+1}^{\lf \frac{r+1}{2}\rf}.
\end{align*}
Moreover, recalling from \cite{DTW23} the definition of the norm in $\mathcal{H}_r^d$:
\begin{align*}
    \Big|\sum_{n=0}^\infty s^n h_n\Big|_{\mathcal{H}_r^d}
=
\inf \left \{ C \ \Big| \ |h_0|, |h_1| \le C, \; \;
|h_n| \le \frac{C\log^d n}{n^{(r+3)/2}}
\quad\text{for all }  n\ge2  \right\}.
\end{align*}
Then from \eqref{eq:locprobdecomp} follows that
\begin{align*}
    \left|\sum_{n=0}^\infty
        \big(
            p_{n}(x)
        -
            p_{n}^{\text{reg}}(x)
        \big) s^n
    \right|_{\mathcal{H}_{r+1}^{\lf \frac{r+1}{2}\rf}}
\le 
    C(1+x)^{r+2}.
\end{align*}
Hence, by \cite[Lemma 4.14]{DTW23}, we have, for {$z \in [1,y]$ and every $j\ge1$},
\begin{align*}
    \left|
    \sum_{k=0}^N
        \big(
            p_{k}(y-z)
        -
            p_{k}^{\text{reg}}(y-z)
        \big)
        a_{n-k}^{({j})}
    \right|
\le
    CA_j \frac{(1+(y-z))^{r+2}}{n^{(r+1+2j)/2}}
    \log^{{\lf \frac{r+1}{2}\rf+\delta_{r+1}}}n,
\end{align*}
where $\delta_m = \mathds{1}_{\{m \text{ is odd} \}}$. Then it holds
\begin{align} \label{eq:convwithanj}
    \left|
    \sum_{k=0}^N
        \big(
            p_{k}(y-z)
        -
            p_{k}^{\text{reg}}(y-z)
        \big)
        a_{n-k}^{(j)}
    \right|
\le
    CA_j \frac{(1+x)^{r+2}}{n^{(r+1+2j)/2}}
    \log^{{\lf \frac{r+2}{2}\rf}}n.
\end{align}
By \eqref{eq:qpower}, $|q^*_{i}(x-z)| \le (1+x)^i$. Using again $[s^n](1-s)^i = 0$ for integer $i$ and $n > i$, we obtain, for $n > r+1$,
\begin{align*}
    |p_{n}^{\text{reg}*}(x-z)| 
\le
    \sum_{\ell=1}^{\lf \frac{r+2}{2}\rf}
        |q^*_{2\ell-1}(x-z)a_n^{(\ell+1)}|\\
\le
    C
    \sum_{\ell=1}^{\lf \frac{r+2}{2}\rf}
        (1+x)^{2\ell-1}|a_n^{(\ell+1)}|.
\end{align*}
Hence, applying \eqref{eq:convwithanj} with $j = \ell+1$ and $x \le \sqrt{n}$, we get
\begin{align*}
    \left|
    \sum_{k=0}^N
        \big(
            p_{k}(y-z)
        -
            p_{k}^{\text{reg}}(y-z)
        \big)
        q^*_{2\ell-1}(x-z)
        a_{n-k}^{(\ell+1)}
    \right|
\le&
    C\frac{(1+x)^{r+2\ell+1}}{n^{(r+2\ell+3)/2}}
    \log^{{\lf \frac{r+2}{2}\rf}}n\\
\le&
    C\frac{(1+x)^{r+2}}{n^{(r+4)/2}}
    \log^{{\lf \frac{{r+2}}{2}\rf}}n.
\end{align*}
Summing over $\ell$ and $z$, we obtain
\begin{align*}
    \left|
    \sum_{z=1}^y
    \sum_{k=0}^N
        \big(
            p_{k}(y-z)
        -
            p_{k}^{\text{reg}}(y-z)
        \big)
        p_{n-k}^{\text{reg}*}(x-z)
    \right|
\le
    C \frac{y(1+x)^{r+2}}{n^{(r+4)/2}}
    \log^{{\lf \frac{r+2}{2}\rf}}n.
\end{align*}
Combining the last inequality and \eqref{eq:lem8firstsum} with \eqref{eq:splitregrest}, we get
\begin{multline*}
    \left|
        \sum_{z=1}^{y}
        \bigg[
        \sum_{k=0}^N
            p_k(y-z)p_{n-k}^*(x-z)
    -
        \sum_{k=0}^N
            p_k^{\text{reg}}(y-z)
            p_{n-k}^{\text{reg}{*}}(x-z)
        \bigg]
    \right|\\
\le
    C\frac{yx^{r+2}}{n^{(r+4)/2}}\log^{\lf \frac{r+2}{2}\rf} n.
\end{multline*}
The sum for $k$ from $N+1$ to $n$ can be done similarly{, with the roles of $p_k(y-z)$ and $p^*_{n-k}(x-z)$ interchanged: 
one uses \eqref{eq:greenbound2} and the bound $(1+y)^{r+2} \le (1+x)^{r+2}$ for the norm of the remainder of $p_k(y-z)$}. 
Thus the lemma is proven.

\end{proof}

\section{Preliminary estimates for the exit-time decomposition}
\label{sec:exitestimates}

In the proofs of Theorems~\ref{thm:2} and \ref{thm:3}, we expand
the unconditional transition probability after the first exit.
Here we estimate the error from discarding large overshoots and
the remainder of the unconditional expansion. We obtain one bound
with a factor $x$ and another independent of the starting point.
Throughout this section, $\E|X_1|^{r+3}<\infty$ for an integer
$r\ge0$, and $\hat y=y/\sigma$, $\hat z=z/\sigma$.

For $x\ge1$ and $y\ge1$,
\begin{align} \label{eq:exitreflection}
    \pr(x+S_n=y,\tau_x>n)
=
    \pr(S_n=y-x)-\pr(x+S_n=y,\tau_x<n).
\end{align}
By the Markov property at $\tau_x$, for every integer $y$,
\begin{align} \label{eq:exitdecomposition}
    \pr(x+S_n=y,\tau_x<n)
=
    \sum_{k=1}^{n-1}\sum_{z\ge0}
        \pr(x+S_k=-z,\tau_x=k)\pr(S_{n-k}=y+z).
\end{align}
We use the expansion from the proof of \cite[Proposition~14]{DTW25a}:
\begin{align} \label{eq:exitunconditional}
    \pr(S_n=y)
=
    e^{-\hat y^2/(2n)}
    \sum_{j=0}^{2r+2}\frac{P_j(\hat y)}{n^{j+1/2}}+w_n(y),
\qquad
    |w_n(y)|\le
    \frac{\varepsilon(n)}{(\sqrt n+|y|)^{r+2}},
\end{align}
where $\varepsilon(n)\downarrow0$, the polynomials $P_j$ have
degree at most $\lf3j/2\rf$, and $P_0=1/(\sigma\sqrt{2\pi})$.
Retaining its polynomial part for $0\le z\le |y|/4$, we define
\begin{align} \label{eq:exitPndef}
    \mathcal P_n(x,y)
:=
    \sum_{k=1}^{n-1}\sum_{0\le z\le |y|/4}
    \pr(x+S_k=-z,\tau_x=k)
    \bigl(\pr(S_{n-k}=y+z)-w_{n-k}(y+z)\bigr).
\end{align}
Thus the error consists of two sums:
\begin{align} \label{eq:exiterroridentity}
&\pr(x+S_n=y,\tau_x<n)-\mathcal P_n(x,y)\nonumber\\
&\hspace{2cm}\quad=
    \sum_{k=1}^{n-1}\sum_{z>|y|/4}
    \pr(x+S_k=-z,\tau_x=k)\pr(S_{n-k}=y+z)\nonumber\\
&\hspace{2cm}\qquad+
    \sum_{k=1}^{n-1}\sum_{0\le z\le |y|/4}
    \pr(x+S_k=-z,\tau_x=k)w_{n-k}(y+z).
\end{align}

To bound these sums, we first estimate the walk before exit.
Recall that $S^*$ is the reflected walk, whose increments have
the law of $-X_1$, and $\tau_u^*$ is its first exit time from
$(0,\infty)$ when started at $u$. Also,
$\overline\tau=\inf\{n\ge1:S_n<0\}$ and $\tau^*=\tau_0^*$.

\begin{lemma} \label{lem:exitwalkbounds}
For $x,u\ge1$ and $n\ge1$,
\begin{align} \label{eq:exitlocalbound}
    \pr(x+S_n=u,\tau_x>n)
    \le \frac{Cu}{n}\min\left(1,\frac{x}{\sqrt n}\right).
\end{align}
Moreover, for $x,u\ge1$,
\begin{align} \label{eq:exitgreenbound}
    \sum_{n=0}^{\infty}\pr(x+S_n=u,\tau_x>n)
    \le C\min(x,u).
\end{align}
\end{lemma}
\begin{proof}
The bound \eqref{eq:exitlocalbound} follows by combining
\eqref{eq:basice_xy_estim} and \eqref{eq:oneendpointbound}.

For \eqref{eq:exitgreenbound}, split the path at the first
time it attains its minimum and reverse the steps before that
time. This gives, for $n\ge0$,
\begin{multline} \label{eq:exitminimum}
    \pr(x+S_n=u,\tau_x>n)\\
=
    \sum_{z=1}^{x\wedge u}\sum_{k=0}^{n}
    \pr(S_k^*=x-z,\tau^*>k)
    \pr(S_{n-k}=u-z,\overline\tau>n-k).
\end{multline}
By the duality lemma \cite[Section~XII.2]{Feller71}, the two
sums below are the renewal masses of the weak ascending ladder
heights of $S$ and the strict ascending ladder heights of $S^*$,
respectively. Strict ladder heights visit each level at most once.
For weak ladder heights, repetitions of a level have a geometric
distribution with a finite mean independent of the level
\cite[Section~XII.1]{Feller71}. Thus
\begin{align} \label{eq:exitrenewalmasses}
    \sum_{k=0}^{\infty}\pr(S_k=a,\overline\tau>k)\le C,
    \qquad
    \sum_{k=0}^{\infty}\pr(S_k^*=a,\tau^*>k)\le C,
    \qquad a\ge0.
\end{align}
Summing \eqref{eq:exitminimum} over $n$, we therefore obtain
\begin{align*}
    \sum_{n=0}^{\infty}\pr(x+S_n=u,\tau_x>n)
&=\sum_{z=1}^{x\wedge u}
    \left(\sum_{k=0}^{\infty}
        \pr(S_k^*=x-z,\tau^*>k)\right)\\
&\qquad\times
    \left(\sum_{j=0}^{\infty}
        \pr(S_j=u-z,\overline\tau>j)\right)
\le C\min(x,u).
\end{align*}
\end{proof}

We now pass from the position before exit to the overshoot.
Conditioning on the last step gives, for $k\ge1$ and $z\ge0$,
\begin{align} \label{eq:exitlaststep}
    \pr(x+S_k=-z,\tau_x=k)
=
    \sum_{u\ge1}\pr(x+S_{k-1}=u,\tau_x>k-1)
        \pr(X_1=-u-z).
\end{align}
For $k=1$, the sum reduces to $u=x$.
Applying \eqref{eq:exitlocalbound} for $k\ge2$ yields
\begin{align} \label{eq:exitpointbound}
    \pr(x+S_k=-z,\tau_x=k)
    \le \frac{Cx}{k^{3/2}}\E[-X_1;-X_1>z].
\end{align}
To bound the full overshoot moment, sum over the exit time and
use \eqref{eq:exitgreenbound} with $\min(x,u)\le x$:
\begin{align} \label{eq:exitmomentbound}
    \E(-x-S_{\tau_x})^{r+2}
&=\sum_{u\ge1}\sum_{k\ge0}
    \pr(x+S_k=u,\tau_x>k)
    \E[(-X_1-u)^{r+2};-X_1>u]\nonumber\\
&\le Cx\sum_{u\ge1}\E[(-X_1-u)^{r+2};-X_1>u]
\le Cx\E|X_1|^{r+3}\le Cx.
\end{align}

\begin{lemma} \label{lem:exitweightedbounds}
Let $g$ be a nonnegative function such that
\begin{align} \label{eq:exitAg}
    A_g:=\sum_{u\ge1}u\sum_{z\ge0}g(z)\pr(X_1=-u-z)<\infty.
\end{align}
Then, uniformly in $x\ge1$,
\begin{align} \label{eq:exitweightedbounds}
    \sum_{k\ge1}\sum_{z\ge0}g(z)
        \pr(x+S_k=-z,\tau_x=k)&\le CA_g,\\
    \sum_{z\ge0}g(z)\pr(x+S_k=-z,\tau_x=k)
        &\le \frac Ck A_g,\qquad k\ge1,
        \label{eq:exitweightedlocalbound}
\end{align}
where $C$ depends only on the distribution of $X_1$.
\end{lemma}
\begin{proof}
Multiplying \eqref{eq:exitlaststep} by $g(z)$, summing over
$z$ and $k$, and using \eqref{eq:exitgreenbound} with
$\min(x,u)\le u$, we obtain
\begin{align*}
&\sum_{k\ge1}\sum_{z\ge0}g(z)
    \pr(x+S_k=-z,\tau_x=k)\\
&\hspace{2cm}\quad=\sum_{u\ge1}
    \left(\sum_{k\ge0}\pr(x+S_k=u,\tau_x>k)\right)
    \sum_{z\ge0}g(z)\pr(X_1=-u-z)\\
&\hspace{2cm}\quad\le C\sum_{u\ge1}u
    \sum_{z\ge0}g(z)\pr(X_1=-u-z)=CA_g.
\end{align*}
For $k\ge2$, \eqref{eq:exitlocalbound} similarly gives
\begin{align*}
    \sum_{z\ge0}g(z)\pr(x+S_k=-z,\tau_x=k)
\le \frac Ck\sum_{u\ge1}u
    \sum_{z\ge0}g(z)\pr(X_1=-u-z)
=\frac Ck A_g.
\end{align*}
For $k=1$, retain only the term $u=x$ in $A_g$ and use $x\ge1$.
\end{proof}

We next apply these bounds to the two sums in
\eqref{eq:exiterroridentity}.
\begin{lemma} \label{lem:exitapproximation}
Uniformly in $x\ge1$ and $y\in\mathbb Z$,
\begin{multline} \label{eq:exitapproximationx}
    \pr(x+S_n=y,\tau_x<n)
=
    \mathcal P_n(x,y)
+o\left(\frac1{(\sqrt n+|y|)^{r+2}}\right)\\
+O\left(\frac{x}{\min(\sqrt n,1+|y|)^{r+3}}\right).
\end{multline}
Moreover, for $n\ge2$,
\begin{align} \label{eq:exitapproximationuniform}
    \left|\pr(x+S_n=y,\tau_x<n)-\mathcal P_n(x,y)\right|
\le \frac C{\sqrt n}\E[X_1^2;X_1<-|y|/4]
    +\frac C{(1+|y|)^{r+2}}.
\end{align}
\end{lemma}
\begin{proof}
We first bound the sum over $z>|y|/4$ in
\eqref{eq:exiterroridentity}. For $k\le n/2$, using the bound $\pr(S_{n} = x) \le C/\sqrt{n}$ and \eqref{eq:exitmomentbound} to obtain
\begin{align*}
\sum_{k=1}^{\lf n/2\rf}\sum_{z>|y|/4}
    \pr(x+S_k=-z,\tau_x=k)\pr(S_{n-k}=y+z)
&\le \frac C{\sqrt n}
    \pr(-x-S_{\tau_x}>|y|/4)\\
&\le 
    \frac{Cx}{\sqrt n(1+|y|)^{r+2}}.
\end{align*}
For $k>n/2$, summing \eqref{eq:exitpointbound} over $z>|y|/4$
and applying the moment assumption gives
\begin{align*}
    \pr(x+S_k<-|y|/4,\tau_x=k)
\le \frac{Cx}{k^{3/2}}\E[X_1^2;-X_1>|y|/4]
\le \frac{Cx}{k^{3/2}(1+|y|)^{r+1}}.
\end{align*}
Using again uniform local bound $\pr(S_{n} = x) \le C/\sqrt{n}$ and summing
over $k>n/2$, we obtain
\begin{multline} \label{eq:exitlargetailx}
    \sum_{k=1}^{n-1}\sum_{z>|y|/4}
    \pr(x+S_k=-z,\tau_x=k)\pr(S_{n-k}=y+z)\\
\le Cx\left(\frac1{\sqrt n(1+|y|)^{r+2}}
            +\frac1{n(1+|y|)^{r+1}}\right)
\le \frac{Cx}{\min(\sqrt n,1+|y|)^{r+3}}.
\end{multline}

For a bound independent of $x$, we use Lemma~\ref{lem:exitweightedbounds}
with
$g(z)=\mathbf1_{\{z>|y|/4\}}$. In this case
\begin{align*}
    A_g=\sum_{u\ge1}u\pr(X_1<-u-|y|/4)
    \le \E[X_1^2;X_1<-|y|/4].
\end{align*}
For $k\le n/2$, use the bound \eqref{eq:exitweightedbounds};
for $k>n/2$, use the bound \eqref{eq:exitweightedlocalbound}, 
together with uniform local bound, these give
\begin{multline} \label{eq:exitlargetailuniform}
    \sum_{k=1}^{n-1}\sum_{z>|y|/4}
    \pr(x+S_k=-z,\tau_x=k)\pr(S_{n-k}=y+z)\\
\le \frac{CA_g}{\sqrt n}
    +\frac{CA_g}{n}\sum_{k=\lf n/2\rf+1}^{n-1}
         \frac1{\sqrt{n-k}}
\le \frac C{\sqrt n}\E[X_1^2;X_1<-|y|/4].
\end{multline}

We next estimate the last sum in \eqref{eq:exiterroridentity}.
For $0\le z\le |y|/4$, we have $|y+z|\ge3|y|/4$, hence for such $z$ and $k \le n/2$ 
we deduce from\eqref{eq:exitunconditional} the bound
\begin{align*}
    |w_{n-k}(y+z)| \le \frac{C\varepsilon(\lceil n/2\rceil)}
          {(\sqrt n+|y|)^{r+2}}
\end{align*}
and for the same $z$ and $k \ge n/2$
\begin{align} \label{eq:wspatialbound}
    |w_{n-k}(y+z)| \le \frac{C}{(1+|y|)^{r+2}}.
\end{align}
Summing \eqref{eq:exitpointbound} over $z\ge0$ one 
obtains $\pr(\tau_x=k)\le Cxk^{-3/2}$ for $k\ge2$.
Using this bound, the bound $\sum_k\pr(\tau_x=k)\le1$ and applying 
the estimations on $w_{n-k}(y+z)$ for the corresponding $k$, we obtain
\begin{align} \label{eq:exitwboundx}
&\left|\sum_{k=1}^{n-1}\sum_{0\le z\le |y|/4}
    \pr(x+S_k=-z,\tau_x=k)w_{n-k}(y+z)\right|\nonumber\\
&\hspace{1.5cm}\quad\le \frac{C\varepsilon(\lceil n/2\rceil)}
          {(\sqrt n+|y|)^{r+2}}
    \sum_{k=1}^{\lf n/2\rf}\pr(\tau_x=k)
    +\frac{C}{(1+|y|)^{r+2}}
    \sum_{k=\lf n/2\rf+1}^{n-1}\pr(\tau_x=k)\nonumber\\
&\hspace{1.5cm}\quad\le \frac{C\varepsilon(\lceil n/2\rceil)}
          {(\sqrt n+|y|)^{r+2}}
    +\frac{Cx}{\sqrt n(1+|y|)^{r+2}}.
\end{align}
Alternatively, bound \eqref{eq:wspatialbound} and
$\sum_k\pr(\tau_x=k)\le1$ give
\begin{align} \label{eq:exitwbounduniform}
    \left|\sum_{k=1}^{n-1}\sum_{0\le z\le |y|/4}
    \pr(x+S_k=-z,\tau_x=k)w_{n-k}(y+z)\right|
    \le \frac C{(1+|y|)^{r+2}}.
\end{align}
Substitute \eqref{eq:exitlargetailx} and \eqref{eq:exitwboundx}
into \eqref{eq:exiterroridentity} to obtain \eqref{eq:exitapproximationx}.
Using \eqref{eq:exitlargetailuniform} and
\eqref{eq:exitwbounduniform} instead proves
\eqref{eq:exitapproximationuniform}.
\end{proof}

\section{Proof of Theorem \ref{thm:2}}\label{sec:3}
We follow the proof of Theorem 1 in \cite{DTW25a}.
The preliminary estimates reduce the problem to expanding
$\mathcal P_n(x,y)$. To expand the overshoot moments for the
walk started at $x$, we use Theorem~\ref{thm:1}.

Throughout this section, $\hat y=y/\sigma$ and $\hat z=z/\sigma$.
By \eqref{eq:exitunconditional}, \eqref{eq:exitPndef}, and
Lemma~\ref{lem:exitapproximation}, uniformly in $x\ge1$ and
$y\in\mathbb Z$,
\begin{align} \label{eq:approxPnxythm2}
    \pr(x+S_n=y,\tau_x<n)
&=\sum_{j=0}^{2r+2}\sum_{k=1}^{n-1}\sum_{0\le z\le |y|/4}
    \pr(x+S_k=-z,\tau_x=k)\nonumber\\
&\qquad{}\times e^{-\frac{(\hat y+\hat z)^2}{2(n-k)}}
    \frac{P_j(\hat y+\hat z)}{(n-k)^{j+1/2}}\nonumber\\
&\quad+o\left(\frac1{(\sqrt n+|y|)^{r+2}}\right)
    +O\left(\frac{x}{\min(\sqrt n,1+|y|)^{r+3}}\right).
\end{align}

We first expand the kernel in powers of the overshoot, retaining
terms up to degree $r+1$. Write
$P_j(t)=\sum_{q=0}^{\lf3j/2\rf}a_{q,j}t^q$. For
$0\le j\le2r+2$ and $0\le s\le q\le\lf3j/2\rf$, set
\begin{align} \label{eq:Tdefthm2}
    T_{j,q,s}(x,y)
:=\hat y^{q-s}\sum_{k=1}^{n-1}\sum_{0\le z\le |y|/4}
    \hat z^s e^{-\frac{(\hat y+\hat z)^2}{2(n-k)}}
    \frac{\pr(x+S_k=-z,\tau_x=k)}{(n-k)^{j+1/2}}.
\end{align}
Expanding each power of $\hat y+\hat z$ in
\eqref{eq:exitPndef} gives
\begin{align} \label{eq:PnTsumthm2}
    \mathcal P_n(x,y)
    =\sum_{j=0}^{2r+2}\sum_{q=0}^{\lf3j/2\rf}\sum_{s=0}^{q}
        a_{q,j}\binom qs T_{j,q,s}(x,y).
\end{align}

We first bound the terms with $s\ge r+2$ using the overshoot's 
moment of order $r+2$.

For $s \ge r+2$ and $0 \le z \le |y|/4$ we have 
$|y|^{q-s}z^s \le C |y|^{q-(r+2)} z^{r+2}$ and $|\hat y + \hat z| \ge 3|\hat y|/4$, hence
\begin{align*}
    |T_{j,q,s}&(x,y)|
\le
    C |y|^{q-(r+2)}
    \sum_{k=1}^{n-1}\sum_{0\le z\le |y|/4}
    z^{r+2} e^{-\frac{(\hat y+\hat z)^2}{2(n-k)}}
    \frac{\pr(x+S_k=-z,\tau_x=k)}{(n-k)^{j+1/2}}\\
&\le 
    C\left(
        \sup\limits_{m\ge 1}
        \frac
        {|y|^{q-(r+2)}e^{-\hat y^2/(8m)}}
        {m^{j+1/2}}
    \right)
    \sum_{k\ge 1}\sum_{0\le z\le |y|/4}
    z^{r+2}
    \pr(x+S_k=-z,\tau_x=k).
\end{align*}
The supremum decay as $(1+|y|)^{-(r+3)}$ since $q \le 2j$. Due to 
\eqref{eq:exitmomentbound} the sum is bounded by $Cx$. Hence, for $s \ge r+2$
\begin{align*}
    |T_{j,q,s}&(x,y)|
\le
    \frac{Cx}{(1+|y|)^{r+3}}.
\end{align*}

The same supremum bound, together with \eqref{eq:exitmomentbound},
gives, for $\alpha\ge0$ and $\alpha-2\beta\le-(r+3)$,
\begin{align} \label{eq:keyboundthm2}
|\hat y|^\alpha\sum_{k=1}^{n-1}
    \frac{e^{-\hat y^2/(8(n-k))}}{(n-k)^\beta}
    \sum_{0\le z\le |y|/4}\hat z^{r+2}
        \pr(x+S_k=-z,\tau_x=k)
\le \frac{Cx}{(1+|y|)^{r+3}},
\end{align}
we need this estimate later.

After bounding terms with $s\ge r+2$ and weakening the reminder, \eqref{eq:PnTsumthm2} becomes
\begin{align} \label{eq:pxydecompthm2}
    \mathcal P_n(x,y)
=\sum_{j=0}^{2r+2}\sum_{q=0}^{\lf3j/2\rf}
    \sum_{s=0}^{q\wedge(r+1)}
    a_{q,j}\binom qs T_{j,q,s}(x,y)
+O\left(\frac{x}{\min(\sqrt n,1+|y|)^{r+3}}\right).
\end{align}

For $s\le q\wedge(r+1)$, factor the Gaussian in $T_{j,q,s}$:
\[
    e^{-\frac{(\hat y+\hat z)^2}{2(n-k)}}
    =e^{-\frac{\hat y^2}{2(n-k)}}
     e^{-\frac{\hat y\hat z}{n-k}}
     e^{-\frac{\hat z^2}{2(n-k)}}.
\]
As in \cite[Section~4]{DTW25a}, expand the last two factors,
indexing their Taylor terms by $\gamma$ and $\mu$,
respectively, and retaining $\gamma+2\mu\le r-s+1$.
On the truncated range, $|\hat y\hat z|\le\hat y^2/4$.
After multiplication by the other factors in $T_{j,q,s}$,
both Taylor remainders are bounded using the moment of
order $r+2$ in \eqref{eq:keyboundthm2}, with
\[
    (\alpha,\beta)=
    (q+r-2s+2,j+r-s+5/2)
    \quad\text{or}\quad
    (q-s+\gamma,j+(r-s+\gamma+3)/2).
\]
In both cases $\alpha\ge0$ and
$\alpha-2\beta=q-2j-r-3\le-(r+3)$.
Thus, uniformly in $x\ge1$ and $y\in\mathbb Z$,
\begin{align} \label{eq:Tjqs:thm2}
    T_{j,q,s}(x,y)
&=\sum_{\gamma+2\mu\le r-s+1}
    \frac{(-1)^{\gamma+\mu}\hat y^{q-s+\gamma}}
         {\gamma!\mu!2^\mu}
    \sum_{k=1}^{n-1}
    \frac{e^{-\hat y^2/(2(n-k))}}
         {(n-k)^{j+\gamma+\mu+1/2}}\nonumber\\
&\hspace{2cm}\qquad{}\times\sum_{0\le z\le |y|/4}
    \hat z^{s+\gamma+2\mu}
    \pr(x+S_k=-z,\tau_x=k)\nonumber\\
&\hspace{5cm}\quad+O\left(\frac{x}{\min(\sqrt n,1+|y|)^{r+3}}\right).
\end{align}
Note that now the overshoot variable $z$ in the formula appears only in monomials.
We extend the inner sum to all $z\ge0$.
On the added tail, $z>|y|/4$ and $z\ge1$, so
$z\ge c(1+|y|)$. Thus, for $0\le h\le r+1$,
\eqref{eq:exitmomentbound} gives
\begin{align} \label{eq:tailsumthm2}
    \sum_{k\ge1}\sum_{z>|y|/4}\hat z^h
        \pr(x+S_k=-z,\tau_x=k)
&\le \frac{C\E(-x-S_{\tau_x})^{r+2}}
             {(1+|y|)^{r-h+2}}
\le \frac{Cx}{(1+|y|)^{r-h+2}}.
\end{align}
For $h=s+\gamma+2\mu\le r+1$, the supremum of the
Gaussian weight in \eqref{eq:Tjqs:thm2} is at most
$C(1+|y|)^{q-2j-h-1}$. Together with
\eqref{eq:tailsumthm2}, this gives
\begin{align*}
&|\hat y|^{q-s+\gamma}\sum_{k=1}^{n-1}
    \frac{e^{-\hat y^2/(2(n-k))}}{(n-k)^{j+\gamma+\mu+1/2}}
    \sum_{z>|y|/4}\hat z^h\pr(x+S_k=-z,\tau_x=k)\\
&\hspace{3cm}\qquad\le \frac{Cx}{(1+|y|)^{r+3+2j-q}}
\le \frac{Cx}{\min(\sqrt n,1+|y|)^{r+3}}.
\end{align*}

Thus we can remove the truncation in \eqref{eq:Tjqs:thm2}
and write, with $h=s+\gamma+2\mu$ in each summand,
\begin{align} \label{eq:Texpectationthm2}
    T_{j,q,s}(x,y)
&=\sum_{\gamma+2\mu\le r-s+1}
    \frac{(-1)^{\gamma+\mu}\hat y^{q-s+\gamma}}
         {\gamma!\mu!2^\mu}\nonumber\\
&\hspace{1cm}\quad{}\times\sum_{k=1}^{n-1}
    \frac{e^{-\hat y^2/(2(n-k))}}
         {(n-k)^{j+\gamma+\mu+1/2}}
    \sigma^{-h}\E[(-x-S_k)^h;\tau_x=k]\nonumber\\
&\hspace{5cm}\quad+O\left(\frac{x}{\min(\sqrt n,1+|y|)^{r+3}}\right).
\end{align}
We next expand the expectations in this formula.
\begin{lemma} \label{lem:momentsoflocprob}
Under the assumptions of Theorem~\ref{thm:2}, for integer
$0\le h\le r+1$ there exist functions $b_\ell^{(h)}(x)$ such that
\begin{align} \label{eq:fullmomentsthm2}
    \sigma^{-h}\E[(-x-S_k)^h;\tau_x=k]
=\sum_{\ell=0}^{\lf(r-h)/2\rf}
    \frac{b_\ell^{(h)}(x)}{k^{\ell+3/2}}+v_k^{(h)}(x),
\end{align}
where, uniformly in $x\ge1$,
\begin{align} \label{eq:momentcoefficientsthm2}
    |b_\ell^{(h)}(x)| \nonumber
    &\le Cx^{2\ell+1}\\
    |v_k^{(h)}(x)|
    &\le \frac{Cx^{r-h+2}\log^{\lf(r-h+2)/2\rf}k}
                   {k^{(r-h+4)/2}}, \quad k\ge2
\end{align}
and additionally $|v_1^{(h)}(x)|\le Cx^{r-h+1}$.
For $h=r+1$, the sum over $\ell$ is empty.
\end{lemma}
\begin{proof}
For $h\le r$, apply Theorem~\ref{thm:1} with $r-h$
in place of $r$. The last-step formula involves $k-1$ steps,
so we expand $(k-1)^{-j-1/2}$ in powers of $k^{-1}$
and denote the resulting coefficients by
$\widetilde U_j$. For $k\ge3$, this gives
\begin{multline} \label{eq:thm1shifted}
    \pr(x+S_{k-1}=u,\tau_x>k-1)
=\sum_{\ell=0}^{\lf(r-h)/2\rf}
    \frac{\widetilde U_{\ell+1}(x,u)}{k^{\ell+3/2}}\\
+O\left(
    \frac{\min(x,u)(x+u)^{r-h+2}\log^{\lf(r-h+2)/2\rf}k}
         {k^{(r-h+4)/2}}
    \right),
\end{multline}
where
$|\widetilde U_j(x,u)|\le C\min(x,u)(x+u)^{2j-1}$.
Multiplying \eqref{eq:thm1shifted} by
$\hat z^h\pr(X_1=-u-z)$ and summing over $u\ge1$ and
$z\ge0$, we obtain \eqref{eq:fullmomentsthm2} with the coefficients
\begin{align} \label{eq:bdefthm2}
    b_\ell^{(h)}(x)
    =\sum_{u\ge1}\sum_{z\ge0}
        \hat z^h\widetilde U_{\ell+1}(x,u)
        \pr(X_1=-u-z).
\end{align}
Using $x+u \le 2 \max(x,u)$ and splitting the sum over $u$ at $u=x$ gives
\begin{align*}
    |b_\ell^{(h)}(x)|
&\le Cx^{2\ell+1}\E|X_1|^{h+2}
    +Cx\E|X_1|^{2\ell+h+3}
\le Cx^{2\ell+1},
\end{align*}
since $2\ell+h\le r$.
For the remainder, the same splitting gives
\begin{align*}
&\sum_{u\ge1}\sum_{z\ge0}\hat z^h
    \min(x,u)(x+u)^{r-h+2}\pr(X_1=-u-z)\\
&\qquad\le Cx^{r-h+2}\E|X_1|^{h+2}
    +Cx\E|X_1|^{r+3}\le Cx^{r-h+2}.
\end{align*}
Together with \eqref{eq:thm1shifted}, this proves the
stated bound for $v_k^{(h)}$.

For $h=r+1$, \eqref{eq:exitpointbound} directly gives
\[
    \sigma^{-(r+1)}\E[(-x-S_k)^{r+1};\tau_x=k]
    \le \frac{Cx}{k^{3/2}}\E|X_1|^{r+3},
    \qquad k\ge2.
\]
For $k=1$,
\[
    \sigma^{-h}\E[(-x-S_1)^h;\tau_x=1]
    \le C\E|X_1|^h\le C.
\]
For $k=2$, \eqref{eq:exitpointbound} gives
\[
    \sigma^{-h}\E[(-x-S_2)^h;\tau_x=2]\le Cx.
\]
Since $2\ell+1\le r-h+1$ for all retained $\ell$,
the coefficient bound gives the stated estimates for
$v_1^{(h)}$ and $v_2^{(h)}$.
\end{proof}

For each $\gamma,\mu$, put $h=s+\gamma+2\mu$.
Substituting \eqref{eq:fullmomentsthm2} into
\eqref{eq:Texpectationthm2} gives
\begin{multline} \label{eq:momentconvolutionsthm2}
    T_{j,q,s}(x,y)
=\sum_{\gamma+2\mu\le r-s+1}
    \frac{(-1)^{\gamma+\mu}\hat y^{q-s+\gamma}}
         {\gamma!\mu!2^\mu}\\
    {}\times\sum_{\ell=0}^{\lf(r-h)/2\rf}b_\ell^{(h)}(x)
    \sum_{k=1}^{n-1}
    \frac{e^{-\hat y^2/(2(n-k))}}
         {k^{\ell+3/2}(n-k)^{j+\gamma+\mu+1/2}}\\
    +R_{j,q,s}(x,y)
    +O\left(\frac{x}{\min(\sqrt n,1+|y|)^{r+3}}\right),
\end{multline}
where
\begin{align} \label{eq:Rdefthm2}
    R_{j,q,s}(x,y)
&=\sum_{\gamma+2\mu\le r-s+1}
    \frac{(-1)^{\gamma+\mu}\hat y^{q-s+\gamma}}
         {\gamma!\mu!2^\mu}
    \sum_{k=1}^{n-1}
    \frac{v_k^{(h)}(x)e^{-\hat y^2/(2(n-k))}}
         {(n-k)^{j+\gamma+\mu+1/2}}.
\end{align}
Next we need the following analoge of \cite[Lemma~20]{DTW25a}, it says 
that $R_{j,q,s}$ can be approximated by powers of $n$ with polynomials in 
$y/\sqrt{n}$ and Gaussian factors.
\begin{lemma} \label{lem:Rexpthm2}
For indices $j,q,s$ in \eqref{eq:pxydecompthm2}, there are polynomials
$H_{\nu,j,q,s}(\cdot,x)$, such that, uniformly in
$x\ge1$ and $y\in\mathbb Z\setminus\{0\}$,
\begin{align} \label{eq:Rexpthm2}
    R_{j,q,s}(x,y)
    =e^{-\frac{\hat y^2}{2n}}\sum_{\nu=1}^{r+2}
        \frac{H_{\nu,j,q,s}(\hat y/\sqrt n,x)}{n^{\nu/2}}
    +O\left(\frac{x^{r+2}\log^{\lf(r+2)/2\rfloor+1}n}
                     {\min(\sqrt n,1+|y|)^{r+3}}\right).
\end{align}
\end{lemma}

\begin{proof}
    It is sufficient to consider each pair $\gamma$ and $\mu$ in \eqref{eq:Rdefthm2}. 
    We claim that for some polynomials $F_\eta(\cdot \, , x)$
\begin{multline} \label{eq:vconvolutionthm2}
    \hat y^{q-s+\gamma}\sum_{k=1}^{n-1}
    \frac{v_k^{(h)}(x)e^{-\hat y^2/(2(n-k))}}
         {(n-k)^{j+\gamma+\mu+1/2}}
=\frac{\hat y^{q-s+\gamma}e^{-\hat y^2/(2n)}}
       {n^{j+\gamma+\mu+1/2}}
    \sum_{\substack{\eta\ge0\\2j-q+h+2\eta\le r+1}}
    \frac{F_\eta(\hat y^2/n,x)}{n^\eta}\\
+O\left(
    \frac{x^{r+2}\log^{\lf(r+2)/2\rf+1}n}
         {\min(\sqrt n,1+|y|)^{r+3}}\right).
\end{multline}
    For $k>\lf n/2\rfloor$, we use
    \eqref{eq:momentcoefficientsthm2} and \cite[Lemma 6]{DTW25a}
    with $a=(r-h+4)/2$ and $b=j+\gamma+\mu+1/2>1$:
\begin{align*}
    \left|
    \hat y^{q-s+\gamma}\sum_{k=\lf n/2 \rf +1 }^{n-1}
    \frac{v_k^{(h)}(x)e^{-\hat y^2/(2(n-k))}}
         {(n-k)^{j+\gamma+\mu+1/2}}
\right|
&\le
    \frac{C x^{r-h+2} |y|^{q-s+\gamma} \log^{\lf \frac{r+2}{2} \rf}n}
    {n^{(r-h+4)/2}|y|^{2j+2\gamma + 2\mu -1}}\\
=
    \frac{C x^{r-h+2} \log^{\lf \frac{r+2}{2} \rf}n}
    {n^{(r-h+4)/2}|y|^{(2j-q) + h - 1}}
&=
    O\left( x^{r-h+2}\frac{\log^{\lf \frac{r+2}{2} \rf}n}{(\min(\sqrt{n},1+|y|))^{r+3}}\right).
\end{align*}
The case $b\le 1$ imply $j=q=s=\gamma=\mu=h=0$ and can be estimated directly
\begin{align*}
    \left|\sum_{k=\lf n/2\rfloor+1}^{n-1}
        \frac{v_k^{(0)}(x)e^{-\hat y^2/(2(n-k))}}{\sqrt{n-k}}\right|
    &\le \frac{Cx^{r+2}\log^{\lf(r+2)/2\rfloor}n}{n^{(r+4)/2}}
        \sum_{m=1}^{\lf n/2\rfloor}m^{-1/2}\\
    &\le Cx^{r+2}\frac{\log^{\lf(r+2)/2\rfloor}n}
        {(\min(\sqrt n,1+|y|))^{r+3}}.
\end{align*}

For $k\le\lf n/2\rfloor$, use the same $a,b$ and put $d=\lf(r-h+2)/2\rfloor$.
The sequence $v_k^{(h)}(x)/x^{r-h+2}$ satisfies the hypotheses of
\cite[Lemma~7]{DTW25a}. Since $a>1$, that lemma gives
\begin{multline*}
    \hat y^{q-s+\gamma}\sum_{k=1}^{\lf n/2\rfloor}
        \frac{v_k^{(h)}(x)e^{-\hat y^2/(2(n-k))}}{(n-k)^b}
    =\frac{\hat y^{q-s+\gamma} e^{-\hat y^2/(2n)}}{n^b}
        \sum_{\eta=0}^{\lf a\rfloor-1}
        \frac{F_\eta(\hat y^2/n,x)}{n^\eta}\\
    +O\left(\frac{x^{r-h+2}\log^{d+\mathbf{1}_{\{a\in\mathbb N\}}}n}
                     {n^{(r+2j-q+3)/2}}\right).
\end{multline*}
Here $F_\eta(\cdot,x)$ has degree at most $\eta$. By the proof of
\cite[Lemma~7]{DTW25a}, its coefficients are finite linear combinations
of the convergent sums $\sum_{k\ge1}k^p v_k^{(h)}(x)$, $p<a-1$,
and are therefore $O(x^{r-h+2})$. For any term with
$2j-q+h+2\eta>r+1$, the Gaussian factor and this coefficient bound give
\[
    \left|\frac{\hat y^{q-s+\gamma} e^{-\hat y^2/(2n)}}{n^{b+\eta}}
        F_\eta(\hat y^2/n,x)\right|
    \le \frac{Cx^{r-h+2}}{n^{(2j-q+h+1+2\eta)/2}}
    \le \frac{Cx^{r-h+2}}{n^{(r+3)/2}}.
\]
Since $q\le2j$, $d+\mathbf{1}_{\{a\in\mathbb N\}}\le
\lf(r+2)/2\rfloor+1$, and $\min(\sqrt n,1+|y|)\le\sqrt n$,
the first half of the sum therefore satisfies
\begin{multline*}
    \hat y^{q-s+\gamma}\sum_{k=1}^{\lf n/2\rfloor}
        \frac{v_k^{(h)}(x)e^{-\hat y^2/(2(n-k))}}{(n-k)^b}
    =\frac{\hat y^{q-s+\gamma} e^{-\hat y^2/(2n)}}{n^b}
        \sum_{\substack{\eta\ge0\\2j-q+h+2\eta\le r+1}}
        \frac{F_\eta(\hat y^2/n,x)}{n^\eta}\\
    +O\left(\frac{x^{r-h+2}\log^{\lf(r+2)/2\rfloor+1}n}
                     {\min(\sqrt n,1+|y|)^{r+3}}\right).
\end{multline*}
Together with the bound for $k>\lf n/2\rfloor$, this proves
\eqref{eq:vconvolutionthm2}. Since
$\hat y^{q-s+\gamma}n^{-j-\gamma-\mu-\eta-1/2}
=n^{-\nu/2}(\hat y/\sqrt n)^{q-s+\gamma}$, where
$\nu=2j-q+h+2\eta+1$, summing \eqref{eq:vconvolutionthm2}
over the finitely many pairs $(\gamma,\mu)$ in \eqref{eq:Rdefthm2}
gives \eqref{eq:Rexpthm2}.
\end{proof}

By \cite[Proposition~9]{DTW25a}, for nonnegative integers $a,b,c$
there are polynomials $P_{\eta,a,b,c}$ of degree $\eta$ such that for some functions
$Q_{a,b,c}$ holds
\begin{align}\label{prop:mainconv}
&\sum_{k=1}^{n-1}
        \frac{\hat y^c}{k^{a+1/2}(n-k)^{b+1/2}}
        e^{-\frac{\hat y^2}{2(n-k)}}
=
    \operatorname{sgn}(y)\sqrt{2\pi}
    e^{-\frac{\hat y^2}{2n}}
    \frac{Q_{a,b,c}(\hat y/\sqrt n)}{n^{a+b-c/2}}
        \\
\nonumber   
&\hspace{1cm}+
    e^{-\frac{\hat y^2}{2n}}
    \frac{\hat y^c}{n^{b+1/2}}
    \sum_{\eta=0}^{\lceil(r+c)/2\rceil-b+1}
        \frac{P_{\eta,a,b,c}(\hat y^2/n)}{n^{\eta}}
+
    O\left(\frac{1}{\min(\sqrt n,1+|y|)^{r+3}}\right).
\end{align}

We apply \eqref{prop:mainconv} to each $k$-sum in the first term of
\eqref{eq:momentconvolutionsthm2} with
$(a,b,c)=(\ell+1,j+\gamma+\mu,q-s+\gamma)$.

Since $2\ell+h\le r$, by \eqref{eq:momentcoefficientsthm2} we have
$|b_\ell^{(h)}(x)|\le Cx^{2\ell+1}\le Cx^{r+1}$.
Thus, after summing over the finitely many indices, the remainder in
\eqref{prop:mainconv} sum up into
$O(x^{r+1}/\min(\sqrt n,1+|y|)^{r+3})$.

By the proposition functions $Q_{a,b,c}$ are polynomials 
that may contain negative powers.
Let $\nu=2(a+b)-c=2+2j-q+h+2\ell$. For $\nu\ge r+3$,
\cite[Lemma~11]{DTW25a} and the coefficient bound give
\begin{align*}
    |b_\ell^{(h)}(x)|e^{-\hat y^2/(2n)}
    \frac{|Q_{a,b,c}(\hat y/\sqrt n)|}{n^{a+b-c/2}}
    \le \frac{Cx^{r+1}}{\min(\sqrt n,1+|y|)^\nu}
    \le \frac{Cx^{r+1}}{\min(\sqrt n,1+|y|)^{r+3}}.
\end{align*}
For the polynomial part of \eqref{prop:mainconv}, write
\[
    \frac{\hat y^c P_{\eta,a,b,c}(\hat y^2/n)}{n^{b+1/2+\eta}}
    =\frac{(\hat y/\sqrt n)^c P_{\eta,a,b,c}(\hat y^2/n)}
            {n^{(1+2b-c+2\eta)/2}}.
\]
Since $2b-c=2j-q+h\ge0$, these terms start at order $n^{-1/2}$.
The Gaussian factor absorbs the polynomial, so terms with
$1+2b-c+2\eta\ge r+3$ are
$O(x^{r+1}n^{-(r+3)/2})$, within the same bound.

Suming \eqref{eq:momentconvolutionsthm2} over $j,q,s$ with a weight $a_{q,j}\binom qs$ 
and using Lemma~\ref{lem:Rexpthm2} to bound $R_{j,q,s}$, we obtain
\begin{multline*}
    \sum_{j=0}^{2r+2}\sum_{q=0}^{\lf3j/2\rf}
    \sum_{s=0}^{q\wedge(r+1)}
    a_{q,j}\binom qs T_{j,q,s}(x,y)
    =\operatorname{sgn}(y)e^{-\hat y^2/(2n)}
    \sum_{\nu=2}^{r+2}
    \frac{Q_\nu(\hat y/\sqrt n,x)}{n^{\nu/2}}\\
    +e^{-\hat y^2/(2n)}
    \sum_{\nu=1}^{r+2}
    \frac{G_\nu(\hat y/\sqrt n,x)}{n^{\nu/2}}
    +O\left(
    \frac{x^{r+2}\log^{\lf(r+2)/2\rfloor+1}n}
         {\min(\sqrt n,1+|y|)^{r+3}}\right).
\end{multline*}
Here the polynomails $G_\nu$ contains the polynomials $P_{\nu, a,b,c}$ 
from \eqref{prop:mainconv} and $H_{\nu, j,q,s}$ from Lemma~\ref{lem:Rexpthm2}. 
The functions $Q_\nu$ for $2\le\nu\le r+2$ are given by
\begin{align} \label{eq:qnu}
    Q_\nu(t,x)
=\sqrt{2\pi}
    \sum_{2+2j-q+h+2\ell=\nu}
    a_{q,j}\binom qs
    \frac{(-1)^{\gamma+\mu}}{2^\mu\gamma!\mu!}
    b_\ell^{(h)}(x)
    Q_{\ell+1,j+\gamma+\mu,q-s+\gamma}(t).
\end{align}
The sum is over nonnegative integers $j,q,s,\gamma,\mu,\ell$ with
$s\le q\le\lf3j/2\rfloor$ and $h=s+\gamma+2\mu$.

By \eqref{eq:pxydecompthm2} and Lemma~\ref{lem:exitapproximation}, this yields
\begin{multline} \label{eq:firstmainresdecomp}
    \pr(x+S_n=y,\tau_x<n)
=\operatorname{sgn}(y)e^{-\hat y^2/(2n)}
    \sum_{\nu=2}^{r+2}
    \frac{Q_\nu(\hat y/\sqrt n,x)}{n^{\nu/2}}\\
+e^{-\hat y^2/(2n)}
    \sum_{\nu=1}^{r+2}
    \frac{G_\nu(\hat y/\sqrt n,x)}{n^{\nu/2}}
+o\left(\frac1{(\sqrt n+|y|)^{r+2}}\right)\\
+O\left(
    \frac{x^{r+2}\log^{\lf(r+2)/2\rf+1}n}
         {\min(\sqrt n,1+|y|)^{r+3}}\right).
\end{multline}

To eliminate $G_\nu$ from the final expansion, we compare
\eqref{eq:firstmainresdecomp} at negative endpoints with
the unconditional probability. For $y<0$, the walk must
have exited by time $n$, so
\begin{align*}
    \pr(x+S_n=y,\tau_x<n)
&=\pr(S_n=y-x)-\pr(x+S_n=y,\tau_x=n),\\
    \pr(x+S_n=y,\tau_x=n)
&\le\frac{Cx}{n^{3/2}(1+|y|)^{r+2}},
\end{align*}
where the bound follows from \eqref{eq:exitpointbound}
and $\E|X_1|^{r+3}<\infty$.

Fix $x\ge1$ and take $y_n<0$ with
$t_n=y_n/(\sigma\sqrt n)\to t<0$. The bound above and
\eqref{eq:firstmainresdecomp} give
\begin{multline*}
    \pr(S_n=y_n-x)
    =e^{-t_n^2/2}\left[
    \sum_{\nu=1}^{r+2}\frac{G_\nu(t_n,x)}{n^{\nu/2}}
    -\sum_{\nu=2}^{r+2}\frac{Q_\nu(t_n,x)}{n^{\nu/2}}
    \right]
    +o(n^{-(r+2)/2}).
\end{multline*}
Applying \cite[Theorem VII.3.16]{Petrov} to $S_n=y-x$ and expanding
in the shift $x/(\sigma\sqrt n)$ gives polynomials $A_\nu(\cdot,x)$ such that,
uniformly in $x\ge1$ and $y\in\mathbb Z$,
\begin{multline*}
    \pr(S_n=y-x)
    =e^{-\hat y^2/(2n)}\sum_{\nu=1}^{r+2}
    \frac{A_\nu(\hat y/\sqrt n,x)}{n^{\nu/2}}\\
    +O\left(\frac{x^{r+2}}{n^{(r+3)/2}}\right)
    +o\left(\frac1{(\sqrt n+|y-x|)^{r+2}}\right).
\end{multline*}
For fixed $x$ and $y=y_n$, both remainders are
$o(n^{-(r+2)/2})$.
Comparing coefficients successively for every $t<0$ yields
\[
    A_1(t,x)=G_1(t,x),\qquad
    A_\nu(t,x)=G_\nu(t,x)-Q_\nu(t,x),
    \quad 2\le\nu\le r+2.
\]
Since $A_\nu(\cdot \, , x)$ and $G_\nu(\cdot \, , x)$ are polynomials, 
$Q_\nu(\cdot \, , x)$ are polynomials as well and the identities hold for all real $t$.

For $1\le x<\sigma\sqrt n$, we have
$\sqrt n+|y-x|\ge c(\sqrt n+|y|)$. Substituting the coefficient
identities into the preceding expansion therefore gives, uniformly
in this range and $y\in\mathbb Z$,
\begin{multline} \label{eq:uncond_x}
    \pr(S_n=y-x)
=e^{-\hat y^2/(2n)}\left[
    \sum_{\nu=1}^{r+2}
    \frac{G_\nu(\hat y/\sqrt n,x)}{n^{\nu/2}}
    -\sum_{\nu=2}^{r+2}
    \frac{Q_\nu(\hat y/\sqrt n,x)}{n^{\nu/2}}
    \right]\\
+O\left(\frac{x^{r+2}}{n^{(r+3)/2}}\right)
+o\left(\frac1{(\sqrt n+|y|)^{r+2}}\right).
\end{multline}

For $1\le x<\sigma\sqrt n$ and $y\ge1$, use
\[
    \pr(x+S_n=y,\tau_x>n)
    =\pr(S_n=y-x)-\pr(x+S_n=y,\tau_x<n).
\]
Substituting \eqref{eq:firstmainresdecomp} and \eqref{eq:uncond_x},
the terms containing $G_\nu$ cancel and with $P_\nu(t,x):=-2Q_\nu(t,x)$, we obtain
\begin{multline} \label{eq:finalexpansionthm2}
    \pr(x+S_n=y,\tau_x>n)
=e^{-\hat y^2/(2n)}
    \sum_{\nu=2}^{r+2}
    \frac{P_\nu(\hat y/\sqrt n,x)}{n^{\nu/2}}\\
+o\left(\frac1{(\sqrt n+y)^{r+2}}\right)
+O\left(
    \frac{x^{r+2}\log^{\lf(r+2)/2\rf+1}n}
         {\min(\sqrt n,y)^{r+3}}\right).
\end{multline}

We a left to verify the degree, parity, and growth asserted in the theorem.
By the explicit formula in \cite[Proposition~9]{DTW25a},
the largest power in $Q_{a,b,c}(t)$ is at most $c+2a-1$,
and every power has parity $c+1$.
For a summand of \eqref{eq:qnu}, the bound $2q\le3j$ gives
\[
    q-s+\gamma+2\ell+1
    \le 3(2j-q+s+\gamma+2\mu+2\ell)+1
    =3\nu-5.
\]
Also, $q-s+\gamma\equiv\nu\pmod2$, so $P_\nu$ contains
$t^i$ only for $i\equiv\nu+1\pmod2$.
Since $2j-q\ge0$, each summand satisfies
$\nu=2+2j-q+h+2\ell\ge2\ell+2$.
The coefficients of $P_\nu$ are therefore finite linear
combinations of $b_\ell^{(h)}(x)$ with $2\ell+1\le\nu-1$.
Consequently, \eqref{eq:momentcoefficientsthm2} yields
\begin{align} \label{eq:Pgrowththm2}
    |P_\nu(t,x)|\le C(1+|t|)^{3\nu-5}x^{\nu-1}.
\end{align}

For $x\ge\sigma\sqrt n$, \eqref{eq:exitLLT} bounds the probability
by $Cn^{-1/2}$. Together with \eqref{eq:Pgrowththm2}, this gives
\begin{align*}
&\pr(x+S_n=y,\tau_x>n)
    +e^{-\hat y^2/(2n)}\sum_{\nu=2}^{r+2}
    \frac{|P_\nu(\hat y/\sqrt n,x)|}{n^{\nu/2}}\\
&\hspace{1cm}\le C\frac{x^{r+2}}{n^{(r+3)/2}}
    \le C\frac{x^{r+2}\log^{\lf(r+2)/2\rfloor+1}n}
             {\min(\sqrt n,y)^{r+3}}.
\end{align*}
Hence \eqref{eq:finalexpansionthm2} holds also in this range.

Finally, taking $\nu=2$ in \eqref{eq:qnu} and using
$a_{0,0}=1/(\sigma\sqrt{2\pi})$ and
$Q_{1,0,0}(t)=-t$,
as computed in \cite[Section~6]{DTW25a}, we obtain
\[
    P_2(t,x)=\frac{2b_0^{(0)}(x)}{\sigma}\,t.
\]
For every fixed $x\ge1$, Lemma~\ref{lem:momentsoflocprob}
with $h=0$ gives
\[
    \pr(\tau_x=k)
    =b_0^{(0)}(x)k^{-3/2}+o(k^{-3/2}).
\]
Consequently,
\[
    \pr(\tau_x>n)
    =\sum_{k>n}\pr(\tau_x=k)
    =\frac{2b_0^{(0)}(x)}{\sqrt n}+o(n^{-1/2}).
\]
Comparing with \eqref{eq:survivalasymptotics} yields
$b_0^{(0)}(x)=V(x)/(\sigma\sqrt{2\pi})$,
which proves the asserted formula for $P_2$ and completes
the proof.

\section{Proof of Theorem \ref{thm:3}}\label{sec:4}
We follow the scheme of \cite[Section 4]{DTW25a} and apply the Markov
property at the time $\tau_x$. The main new ingredient is the
expansion of the overshoot moments obtained from Theorem~\ref{thm:2}
for the time-reversed walk. We first derive the expansion for
$x \ge 1$ and $y \ne 0$; negative $y$ will be used in the last step
to exclude poles in the coefficients.

Throughout this section, $\hat x=x/\sigma$, $\hat y=y/\sigma$,
and $\hat z=z/\sigma$. 
Put $M:=\min(\sqrt n,x,|y|)$, and let
$\mathcal P_n(x,y)$ be defined by \eqref{eq:exitPndef}.
By Lemma~\ref{lem:exitapproximation} and the moment assumption,
uniformly in $x\ge1$ and $y\ne0$,
\begin{align*}
    \left|\pr(x+S_n=y,\tau_x<n)-\mathcal P_n(x,y)\right|
&\le \frac C{\sqrt n}\E[X_1^2;X_1<-|y|/4]
    +\frac C{(1+|y|)^{r+2}}\\
&\le \frac C{\sqrt n\,|y|^{r+1}}+\frac C{|y|^{r+2}}
\le \frac{C\sqrt n}{M^{r+3}}.
\end{align*}
Thus, uniformly in
$x \ge 1$ and $y \ne 0$,
\begin{align} \label{eq:approxPnxy}
    \pr(x+S_n = y, \tau_x < n)
=
    \mathcal{P}_n(x,y)
+
    O\left(\frac{\sqrt{n}}{M^{r+3}}\right).
\end{align}

Recall the notation
\begin{align*}
        T_{j,q,s}(x,y)
:=\hat y^{q-s}\sum_{k=1}^{n-1}\sum_{0\le z\le |y|/4}
    \hat z^s e^{-\frac{(\hat y+\hat z)^2}{2(n-k)}}
    \frac{\pr(x+S_k=-z,\tau_x=k)}{(n-k)^{j+1/2}}.
\end{align*}
and \eqref{eq:PnTsumthm2}
\[
    \mathcal P_n(x,y)
    =\sum_{j=0}^{2r+2}\sum_{q=0}^{\lf3j/2\rfloor}
      \sum_{s=0}^{q}a_{q,j}\binom qs T_{j,q,s}(x,y),
\]
where $a_{q,j}$ are coefficients of $P_j$ from \eqref{eq:exitunconditional}.

To discard the terms with $s\ge r+1$ and bound the Taylor
remainders below, we use the truncated overshoot moments
\begin{align*}
    \Theta_k^{(h)}(x,y)
:=
    \sum_{0 \le z \le |y|/4}
        \hat{z}^h \pr(x+S_k = -z, \tau_x = k).
\end{align*}

Apply Lemma~\ref{lem:exitweightedbounds} with
$g(z)=\hat z^{r+1}\mathbf 1_{\{z\le |y|/4\}}$.
Since $A_g\le C\E|X_1|^{r+3}$, it gives
$\sum_{k\ge1}\Theta_k^{(r+1)}(x,y)\le C$.
For every $m>0$,
$(\hat y^2/m)^\beta e^{-\hat y^2/(8m)}\le C_\beta$.
Using the full sum of the overshoot moments, we obtain
\begin{align*}
    |\hat y|^\alpha\sum_{k=1}^{n-1}
    \frac{e^{-\frac{\hat y^2}{8(n-k)}}}{(n-k)^\beta}
    \Theta_k^{(r+1)}(x,y)
    &\le
    \left(\sup_{m\ge1}
    \frac{|\hat y|^\alpha e^{-\frac{\hat y^2}{8m}}}{m^\beta}\right)
    \sum_{k\ge1}\Theta_k^{(r+1)}(x,y) \le C|\hat y|^{\alpha-2\beta}
\end{align*}
Hence, for
$\alpha \ge 0$ and $\beta > 1$ with $\alpha - 2\beta \le -(r+2)$,
uniformly in $x \ge 1$ and $y \ne 0$,
\begin{align} \label{eq:keybound}
    |\hat{y}|^{\alpha}
\sum_{k=1}^{n-1}
    \frac{e^{-\frac{\hat{y}^2}{8(n-k)}}}{(n-k)^{\beta}}
    \Theta_k^{(r+1)}(x,y)
\le
    \frac{C\sqrt{n}}{M^{r+3}}.
\end{align}
The last inequality uses $M\le |y|$ and $M\le\sqrt n$.
The estimate also holds with $4$ or $2$ in place of $8$ in the
exponent.

Now let $s \ge r+1$. For $z \le |y|/4$ we have
\begin{align*}
    |\hat{z}+\hat{y}| \ge 3|\hat{y}|/4
    \quad \text{and} \quad
    \hat{z}^s \le (|\hat{y}|/4)^{s-r-1}\hat{z}^{r+1},
\end{align*}
so that by \eqref{eq:keybound} with $\alpha = q-r-1$ and $\beta = j+1/2$
\begin{align*}
    |T_{j,q,s}(x,y)|
\le
    C|\hat{y}|^{q-r-1}
    \sum_{k=1}^{n-1}
        \frac{e^{-\frac{\hat{y}^2}{8(n-k)}}}{(n-k)^{j+1/2}}
        \Theta_k^{(r+1)}(x,y)
\le
    \frac{C\sqrt{n}}{M^{r+3}}.
\end{align*}
The parameters $\alpha$ and $\beta$ satisfy its assumptions since $q \ge s \ge r+1$
and $q \le 2j$.
Hence
\begin{align} \label{eq:pxydecomp1}
    \mathcal{P}_n(x,y)
=
    \sum_{j=0}^{2r+2}
    \sum_{q=0}^{\lf \frac{3j}{2}\rf}
    \sum_{s=0}^{q \wedge r}
        a_{q,j}
        \binom{q}{s}
        T_{j,q,s}(x,y)
+
    O\left(\frac{\sqrt{n}}{M^{r+3}}\right).
\end{align}
Next, we expand the exponent in $T_{j,q,s}(x,y)$, $s \le q \wedge r$,
writing
\begin{align*}
    e^{-\frac{(\hat{z}+\hat{y})^2}{2(n-k)}}
=
    e^{-\frac{\hat{y}^2}{2(n-k)}}
    e^{-\frac{\hat{z}\hat{y}}{n-k}}
    e^{-\frac{\hat{z}^2}{2(n-k)}}.
\end{align*}
By the Taylor formula, for $z \le |y|/4$,
\begin{align*}
    e^{-\frac{\hat{z}\hat{y}}{n-k}}
=
        \sum_{\gamma=0}^{r-s}
            \frac{(-1)^\gamma}{\gamma!}
            \frac{\hat{z}^\gamma \hat{y}^\gamma}{(n-k)^\gamma}
+
    O\left(
        e^{\frac{\hat{y}^2}{4(n-k)}}
        \frac{\hat{z}^{r-s+1}|\hat{y}|^{r-s+1}}{(n-k)^{r-s+1}}
    \right),
\end{align*}
where, since $\theta \in (0,1)$, we used the bound
$$e^{-\theta\frac{\hat{z}\hat{y}}{n-k}} \le e^{\frac{\hat{z}|\hat{y}|}{n-k}}
\le 
e^{\frac{\hat{y}^2}{4(n-k)}}.$$ 
We stop at $\gamma = r-s$, so after multiplication by a factor $\hat{z}^s$ 
from the definition of $T_{j,q,s}$, the remainder contains $\hat{z}^{r+1}$. Its
contribution to $T_{j,q,s}(x,y)$ is, as required, at most
$C\sqrt{n}M^{-(r+3)}$ by \eqref{eq:keybound} with
$\alpha = q+r-2s+1$ and $\beta = j+r-s+3/2$. The assumptions hold since
$s \le q \le 2j$ and $s \le r$.

Further, for every integer $R \ge 1$,
\begin{align*}
    e^{-t^2/2}
=
    \sum_{2\mu \le R-1}
        \frac{(-1)^\mu t^{2\mu}}{2^\mu \mu!}
+
    O\left(t^R\right),
    \qquad t \ge 0,
\end{align*}
by the Taylor formula for $t \le 1$ and trivially for $t > 1$. In the
term containing $\hat{z}^{\gamma}$ from the first expansion, we apply
this with $t = \hat{z}/\sqrt{n-k}$ and $R = r-s-\gamma+1$. The
remainder again contains $\hat{z}^{r+1}$, and its contribution is at
most $C\sqrt{n}M^{-(r+3)}$ by \eqref{eq:keybound} with
$\alpha = q-s+\gamma$ and $\beta = (r+2j+\gamma-s+2)/2$. The assumptions
follow from $s \le q \le 2j$ and $r \ge 1$.
Thus both Taylor remainders contribute
$O(\sqrt{n}\,M^{-(r+3)})$. Substituting the truncated expansions, we
obtain
\begin{multline} \label{eq:Tjqs:thm3}
    T_{j,q,s}(x,y)
=
    \sum_{k=1}^{n-1}
    e^{-\frac{\hat{y}^2}{2(n-k)}}
        \sum_{\gamma + 2\mu \le r-s}
            \frac{(-1)^{\gamma+\mu} \hat{y}^{q-s+\gamma}}
            {\gamma!\, \mu!\, 2^\mu\, (n-k)^{j+\gamma+\mu+1/2}}
            \Theta_k^{(s+\gamma+2\mu)}(x,y)\\
+
    O\left(\frac{\sqrt{n}}{M^{r+3}}\right).
\end{multline}
Now the whole appearance of the overshoot $z$ is hidden inside
$\Theta_k^{(h)}(x,y)$ with index $h$ not exceeding $r$, which ensures its finiteness.
So we can proceed with the expansion of $\Theta_k^{(h)}(x,y)$.

On $\{\tau_x = k\}$ the walk is at some $u \ge 1$ at time $k-1$ and
jumps below zero at time $k$. Reversing time on the first $k-1$
steps turns $\{x+S_{k-1} = u, \tau_x > k-1\}$ into the event that
the reflected walk started at $u$ stays positive and ends at $x$,
to which Theorem~\ref{thm:2} applies. A moment of order $h$ of the overshoot
costs $h$ moments of $X_1$, so Theorem~\ref{thm:2} has to be
applied with $r-h$ in place of $r$. The
expansion holds for all $k \ge 2$, but its remainder is small only
for large $k$. Recall
\begin{align*}
    \Theta_k^{(h)}(x,y)
=
    \sum_{0 \le z \le |y|/4}
        \hat{z}^h \pr(x+S_k = -z, \tau_x = k).
\end{align*}
\begin{lemma}\label{lem:momentsoflocprob2}
Let $0 \le h \le r$. There exist polynomials
$F^{(h)}_2, F^{(h)}_3, \ldots, F^{(h)}_{r-h+2}$ with
$\deg F^{(h)}_m \le 3m-5$, containing $t^i$ only for
$i \equiv m+1 \pmod 2$, such that, uniformly in $x \ge 1$,
$y \ne 0$, and $k \ge 2$,
\begin{multline*}
    \Theta_k^{(h)}(x,y)
=
    e^{-\frac{\hat{x}^2}{2k}}
    \sum_{m=2}^{r-h+2}
        \frac{F^{(h)}_m(\hat{x}/\sqrt{k})}{k^{m/2}}\\
+
    o\left(\frac{1}{(\sqrt{k}+x)^{r-h+2}}\right)
+
    O\left(\frac{\log^{\lf \frac{r-h+2}{2} \rf+1} k}
                {\min(\sqrt{k},x,|y|)^{r-h+3}}\right).
\end{multline*}
\end{lemma}
\begin{proof}
By the total probability formula and the time reversal,
\begin{align} \label{eq:thetareversal}
    \pr(x+S_k = -z, \tau_x = k)
=
    \sum_{u=1}^\infty
        \pr(u + S^*_{k-1} = x, \tau^*_u > k-1)
        \pr(X_1 = -u-z),
\end{align}
where $S^*_n \overset{d}{=} -S_n$ and $\tau^*_u$ is the exit time of $u+S^*_n$
from $(0,\infty)$. Theorem~\ref{thm:2} applied to $S^*$ with $r-h$
in place of $r$ gives, uniformly in $u, x \ge 1$ and $k \ge 2$,
\begin{align} \label{eq:thm2reversed}
    \pr(u + S^*_{k-1} = x, \tau^*_u > k-1)
=
    e^{-\frac{\hat{x}^2}{2(k-1)}}
    \sum_{\nu=2}^{r-h+2}
        \frac{P^*_\nu(\hat{x}/\sqrt{k-1},u)}{(k-1)^{\nu/2}}
+
    \rho_k(u,x),
\end{align}
where $P^*_\nu(\cdot,u)$ are polynomials of degree at most $3\nu-5$
containing $t^i$ only for $i \equiv \nu+1 \pmod 2$, with
coefficients bounded by $Cu^{\nu-1}$, and
\begin{align*}
    \rho_k(u,x)
=
    o\left(\frac{1}{(\sqrt{k}+x)^{r-h+2}}\right)
+
    O\left(\frac{u^{r-h+2}\log^{\lf \frac{r-h+2}{2} \rf+1} k}{\min(\sqrt{k},x)^{r-h+3}}\right).
\end{align*}

We re-expand each term in \eqref{eq:thm2reversed} in powers of $1/k$
and stop at order $k^{-(r-h+3)/2}$. This gives
\begin{align*}
    \pr(u + S^*_{k-1} = x, \tau^*_u > k-1)
=
    e^{-\frac{\hat{x}^2}{2k}}
    \sum_{m=2}^{r-h+2}
        \frac{\widetilde F_m(\hat{x}/\sqrt{k},u)}{k^{m/2}}+
    \widetilde\rho_k(u,x).
\end{align*}
The remainder $\widetilde\rho_k(u,x)$ satisfies the same estimate
as $\rho_k(u,x)$. Indeed, the Gaussian factor absorbs the
polynomial growth in $\hat{x}/\sqrt{k}$, so the re-expansion adds
an error $O(u^{r-h+1}k^{-(r-h+3)/2})$, which is absorbed into the
$O$-term of that estimate. Here
$\widetilde F_m(\cdot,u)$ is a polynomial of degree at most $3m-5$
with coefficients bounded by $Cu^{m-1}$. It contains $t^i$ only for
$i \equiv m+1 \pmod 2$, since the re-expansion produces only the
factors $1/k$ by monomials and $t^2/k$ by the exponent.

We now plug this into \eqref{eq:thetareversal}, multiply by
$\hat{z}^h$ and sum over $z \ge 0$ and $u \ge 1$. Since
$u^iz^h \le (u+z)^{i+h}$ and there are at most $u+z$ pairs $(u,z)$
with a given sum, for $0 \le i \le r-h+2$ we have
\begin{align} \label{eq:uzmomentsum}
    \sum_{z=0}^\infty
    \sum_{u=1}^\infty
        \hat{z}^h u^i \pr(X_1 = -u-z)
\le
    C\E|X_1|^{i+h+1}
\le
    C\E|X_1|^{r+3}.
\end{align}
Therefore
\begin{align*}
    F^{(h)}_m(t)
:=
    \sum_{z=0}^\infty
    \sum_{u=1}^\infty
            \hat{z}^h \widetilde F_m(t,u) \pr(X_1 = -u-z)
\end{align*}
is a polynomial of degree at most $3m-5$. Summing the remainder
with the same weights and using \eqref{eq:uzmomentsum}, we obtain
the expansion of the full overshoot moment:
\begin{multline*}
    \sum_{z=0}^\infty
        \hat{z}^h\pr(x+S_k=-z,\tau_x=k)
=
    e^{-\frac{\hat{x}^2}{2k}}
    \sum_{m=2}^{r-h+2}
        \frac{F_m^{(h)}(\hat{x}/\sqrt{k})}{k^{m/2}}\\
+
    o\left(\frac{1}{(\sqrt{k}+x)^{r-h+2}}\right)
+
    O\left(\frac{\log^{\lf \frac{r-h+2}{2} \rf+1} k}{\min(\sqrt{k},x)^{r-h+3}}\right).
\end{multline*}
To recover $\Theta_k^{(h)}(x,y)$, we remove the terms with
$z>|y|/4$. By \eqref{eq:exitweightedlocalbound} with
$g(z) = \hat{z}^h\mathbf{1}_{\{z > |y|/4\}}$, for which
\begin{align*}
    A_g
\le
    C\E\left[|X_1|^{h+2}; |X_1| > |y|/4\right]
\le
    \frac{C\E|X_1|^{r+3}}{|y|^{r-h+1}},
\end{align*}
the omitted tail satisfies
\begin{align*}
    \sum_{z>|y|/4}\hat z^h\pr(x+S_k=-z,\tau_x=k)
    \le \frac{C}{k|y|^{r-h+1}}
    \le \frac{C}{\min(\sqrt k,x,|y|)^{r-h+3}}.
\end{align*}
This tail is absorbed into the
$O$-term stated in the lemma, which completes the proof.
\end{proof}

For $k \ll x^2$, the remainder of Lemma~\ref{lem:momentsoflocprob2}
is not smaller than its main terms, so the lemma does not give the
required estimate in this range.
Instead, we use that the walk started at $x$ is unlikely to exit
the positive half-line before the time $x^2/\log n$.

\begin{lemma}\label{lem:smallk}
Let $0\le h\le r$ and
\begin{align*}
    k_0 := \min(\lceil cx^2/\log n \rceil, \lf n/2 \rf),
\end{align*}
where $c>0$ is sufficiently small. Then, uniformly in $x\ge1$,
$y\ne0$ and $n\ge3$,
\begin{align} \label{eq:smallkthm3}
    \sum_{k=1}^{k_0}
        \Theta_k^{(h)}(x,y)
\le
    \frac{C}{x^{r+1-h}}.
\end{align}
\end{lemma}
\begin{proof}
If $k_0 = 1$, the left-hand side is
$$\Theta_1^{(h)}(x,y) \le \sigma^{-h}\E[|X_1|^h; X_1 \le -x]
\le Cx^{-(r+3-h)}.$$
Let $k_0 \ge 2$; then $k_0 \le 2cx^2/\log n$. On $\{\tau_x = k\}$
we have $0 \le -x-S_k \le |X_k|$, and therefore
\begin{align*}
    \sum_{k=1}^{k_0}
        \Theta_k^{(h)}(x,y)
&\le
    \sigma^{-h}
    \sum_{k=1}^{k_0}
        \E\big[|X_k|^h; \tau_x = k\big]\\
&\le
    \sigma^{-h}
    \Big(
        x^h\,\pr(\tau_x \le k_0)
    +
        k_0\,\E\big[|X_1|^h; |X_1| > x\big]
    \Big).
\end{align*}
Since $k_0 \le C x^2$, the second term is at most
$k_0x^{-(r+3-h)}\E|X_1|^{r+3} \le Cx^{-(r+1-h)}$.
By Etemadi inequality (see, for example, \cite[Theorem 22.5]{Billingsley1995})
\begin{align*}
    \pr(\tau_x \le k_0) 
\le
    \pr(\max_{k \le k_0} |S_k| \ge x)
\le
    3 \max_{k \le k_0} \pr(|S_k| \ge x/3).
\end{align*}
Then by the Fuk--Nagaev inequality \cite[Corollary 1.8]{Nagaev79}
gives, with some $c_0>0$,
\begin{align*}
    \pr(\tau_x \le k_0)
    \le
    \frac{Ck_0\,\E|X_1|^{r+3}}{x^{r+3}}
    + 6e^{-c_0x^2/k_0}.
\end{align*}
The polynomial term is at most $Cx^{-(r+1)}$ since $k_0\le Cx^2$.
For the exponential term, we apply $x^2/k_0\ge\log n/(2c)$ for first half of the exponent and
$k_0\le n/2$ for the second, this gives, for $c_0/(4c)\ge(r+1)/2$,

\begin{align*}
    e^{-c_0x^2/k_0}
=
    e^{-c_0x^2/2k_0}
    e^{-c_0x^2/2k_0}
\le
    n^{-c_0/(4c)}e^{-c_0x^2/n}
\le 
    \frac{C}{x^{r+1}},
\end{align*}
where the last inequality follows from
$x^{r+1}e^{-c_0x^2/n}\le Cn^{(r+1)/2}$. Thus
\begin{align*}
    x^h\pr(\tau_x \le k_0) \le \frac{C}{x^{r+1-h}},
\end{align*}
which finishes the proof.
\end{proof}

We return to \eqref{eq:Tjqs:thm3}. Fix admissible indices
$j,q,s,\gamma,\mu$, and write $h=s+\gamma+2\mu$ and
$b=j+\gamma+\mu+1/2$. Up to the factor
$(-1)^{\gamma+\mu}(\operatorname{sgn}y)^{q-s+\gamma}/
(\gamma!\,\mu!\,2^\mu)$, whose absolute value is at most $1$,
the weight of $\Theta_k^{(h)}(x,y)$ is
\begin{align*}
    \lambda_k
    :=\frac{e^{-\frac{\hat y^2}{2(n-k)}}
                  |\hat y|^{q-s+\gamma}}
             {(n-k)^b},
    \qquad 1\le k<n.
\end{align*}
We suppress the dependence on $n$, $y$, and the fixed indices.
For $t=|\hat y|/\sqrt{n-k}$, since $q\le2j$ and
$t^{q-s+\gamma+h}e^{-t^2/2}$ is bounded,
\begin{align*}
    \sqrt{n-k}\,|y|^h\lambda_k
    =\sigma^h(n-k)^{-(2j-q)/2}
      t^{q-s+\gamma+h}e^{-t^2/2}
    \le C.
\end{align*}
As $M\le|y|$, this gives
\begin{align} \label{eq:weightthm3}
    \lambda_k\le\frac{C}{\sqrt{n-k}\,M^h},
    \qquad 1\le k<n.
\end{align}
Since $k_0\le n/2$, Lemma~\ref{lem:smallk} therefore yields
\begin{align} \label{eq:smallkoriginalthm3}
    \sum_{k=1}^{k_0} \lambda_k\Theta_k^{(h)}(x,y)
\le
    \frac{C}{\sqrt n\,M^h x^{r+1-h}}
\le
    \frac{C}{\sqrt{n}\,M^{r+1}}
\le
    \frac{C\sqrt{n}}{M^{r+3}},
\end{align}
since $M\le\min(x,\sqrt n)$.

Let $R_k^{(h)}(x,y)$ denote the remainder in
Lemma~\ref{lem:momentsoflocprob2}, that is,
\begin{align*}
    R_k^{(h)}(x,y)
    := \Theta_k^{(h)}(x,y)
    - e^{-\frac{\hat{x}^2}{2k}}
      \sum_{m=2}^{r-h+2}
          \frac{F_m^{(h)}(\hat{x}/\sqrt{k})}{k^{m/2}},
    \qquad k\ge2.
\end{align*}
By \eqref{eq:smallkoriginalthm3}, we can discard $k\le k_0$
in \eqref{eq:Tjqs:thm3}. Applying Lemma~\ref{lem:momentsoflocprob2} for
$k>k_0$, we obtain
\begin{multline} \label{eq:Tjqscutthm3}
    T_{j,q,s}(x,y)
    =\sum_{\gamma+2\mu\le r-s}
      \frac{(-1)^{\gamma+\mu}\hat y^{q-s+\gamma}}
           {\gamma!\,\mu!\,2^\mu}
      \sum_{k=k_0+1}^{n-1}
      \frac{e^{-\frac{\hat y^2}{2(n-k)}}}{(n-k)^b}\\
    {}\times\left[
      e^{-\frac{\hat x^2}{2k}}
      \sum_{m=2}^{r-h+2}
      \frac{F_m^{(h)}(\hat x/\sqrt k)}{k^{m/2}}
      +R_k^{(h)}(x,y)\right]
    +O\left(\frac{\sqrt n}{M^{r+3}}\right).
\end{multline}
Here $h=s+\gamma+2\mu$ and $b=j+\gamma+\mu+1/2$ in each
summand. It remains to estimate the sum of the remainders and
restore the main terms with $k\le k_0$.

We first estimate the contribution of $R_k^{(h)}(x,y)$.
For $k>k_0$, either
$k>cx^2/\log n$ or $k>n/2$. Hence, for $n\ge3$,
\begin{align*}
    \min(\sqrt k,x,|y|)
\ge
    \frac{\min(\sqrt n,x,|y|)}{C\sqrt{\log n}}
=
    \frac{M}{C\sqrt{\log n}}.
\end{align*}
Lemma~\ref{lem:momentsoflocprob2} and $\sqrt k+x\ge M$ therefore
give, for some fixed $D = D(r)$,
\begin{align*}
    |R_k^{(h)}(x,y)|
    \le \frac{\delta(k)}{M^{r-h+2}}
       +\frac{C\log^D n}{M^{r-h+3}},
    \qquad k_0<k<n,
\end{align*}
where $\delta(k)\ge0$ tends to zero independently of $x,y$.
Using
\begin{align*}
    \sum_{k=2}^{n-1}\frac{\delta(k)}{\sqrt{n-k}}=o(\sqrt n),
\end{align*}
\eqref{eq:weightthm3} and $M\le\sqrt n$, we conclude, uniformly in $x\ge1$ and
$y\ne0$,
\begin{align} \label{eq:remaindersumthm3}
\nonumber
    \sum_{k=k_0+1}^{n-1}\lambda_k|R_k^{(h)}(x,y)|
    &\le \frac{C}{M^{r+2}}
          \sum_{k=2}^{n-1}\frac{\delta(k)}{\sqrt{n-k}}
       +\frac{C\sqrt n\log^D n}{M^{r+3}}\\
    &=o\left(\frac{n}{M^{r+3}}\right).
\end{align}

We now restore the main terms with $k\le k_0$ in
\eqref{eq:Tjqscutthm3}. In this range \eqref{eq:weightthm3}
and $M\ge1$ give $\lambda_k\le Cn^{-1/2}$.
The term with $k=1$ is bounded by
$C/(\sqrt n\,x^{r+1})\le C\sqrt n/M^{r+3}$ due to the factor $e^{-\hat{x}^2/2}$.

For $2\le k\le k_0$, we have $x^2/k\ge\log n/(2c)$.
Thus the Gaussian factor absorbs the polynomial and gives
$$e^{-\hat{x}^2/2k}|F^{(h)}_m(\hat{x}/\sqrt{k})|
\le Cn^{-1/(8\sigma^2c)}.$$
Using $\lambda_k\le Cn^{-1/2}$, the sum over $2\le k\le k_0$
is at most $Cn^{1/2-1/(8\sigma^2c)}$. Adding the term $k=1$,
we obtain, for $c$ small enough,
\begin{align} \label{eq:smallkmainthm3}
    \sum_{k=1}^{k_0}
        \lambda_k e^{-\hat{x}^2/2k}
        \frac{|F^{(h)}_m(\hat{x}/\sqrt{k})|}{k^{m/2}}
    \le \frac{C\sqrt{n}}{M^{r+3}}.
\end{align}
By \eqref{eq:remaindersumthm3} and
\eqref{eq:smallkmainthm3}, discarding the remainders and restoring
the main terms for $k\le k_0$ introduces an error $o(n/M^{r+3})$.
The factors omitted from $\lambda_k$ are bounded, and summing
over the finitely many pairs $(\gamma,\mu)$ preserves this bound.
Thus \eqref{eq:Tjqscutthm3} becomes
\begin{multline} \label{eq:Tjqs4:thm3}
    T_{j,q,s}(x,y)
=
        \sum_{\gamma + 2\mu \le r-s}
        \sum_{m=2}^{r-h+2}
            \frac{(-1)^{\gamma+\mu}\hat{y}^{q-s+\gamma}}{\gamma!\, \mu!\, 2^\mu}
    \sum_{k=1}^{n-1}
        \frac{e^{-\frac{\hat{x}^2}{2k} - \frac{\hat{y}^2}{2(n-k)}}
              F^{(h)}_m(\hat{x}/\sqrt{k})}
             {k^{m/2}(n-k)^{j+\gamma+\mu+1/2}}\\
+
    o\left(\frac{n}{M^{r+3}}\right).
\end{multline}

We first discard terms in \eqref{eq:Tjqs4:thm3} that can be
absorbed into the remainder, then replace the remaining sums
over $k$ by integrals.

We control the error through a bound on
the derivatives of the summand. Let
\begin{align*}
    f(t)
:=
    \frac{e^{-\frac{\hat{x}^2}{2t} - \frac{\hat{y}^2}{2(n-t)}}
          F^{(h)}_m(\hat{x}/\sqrt{t})}
         {t^{m/2}(n-t)^{b}},
    \qquad 0 < t < n,
\end{align*}
be the summand of \eqref{eq:Tjqs4:thm3} with $k$ replaced by the
continuous variable $t$. Differentiating term by term and
absorbing the polynomial factors into the Gaussians, we obtain,
for every integer $p\ge0$,
\begin{align} \label{eq:fderthm3}
    |f^{(p)}(t)|
&\le
    \frac{Ce^{-\frac{\hat{y}^2}{4n}}}{n^{b}}
    \frac{e^{-\frac{\hat{x}^2}{4t}}}{t^{p+m/2}},
    \qquad t \le n/2,\\
    |f^{(p)}(t)|
&\le \nonumber
    \frac{C}{n^{m/2}}\,
    \frac{e^{-\frac{\hat{y}^2}{4(n-t)}}}
         {(n-t)^{p+b}},
    \qquad t \ge n/2.
\end{align}

Put $\eta=2j-q+m+h-1$, then with $b = j + \gamma + \mu+1/2$ and $h = s + \gamma + 2\mu$
\begin{align} \label{eq:etarelation}
    \eta + (q-s+\gamma) = m + 2b-2.
\end{align}
We omit the terms in
\eqref{eq:Tjqs4:thm3} with $\eta\ge r+2$.
For these terms, $2\le m\le r-h+2\le\eta$.
Applying \eqref{eq:fderthm3} with $p=0$ and absorbing powers into
the Gaussian factors before summing, we obtain, uniformly in
$x\ge1, y\neq 0$,
\begin{align*}
    |\hat{y}|^{q-s+\gamma}\sum_{k=1}^{n-1}|f(k)|
    &\le \frac{C}{n^{(\eta-m)/2}x^m}
    +\frac{C}{n^{(m-2)/2}|y|^{\eta-m+2}}\\
    &\le \frac{C}{M^\eta}
    \le \frac{C}{M^{r+2}}
    =o\left(\frac{n}{M^{r+3}}\right),
\end{align*}
where we used $M\le\sqrt n,x,|y|$ and \eqref{eq:etarelation}.

It remains to approximate the sums with $\eta\le r+1$ by integrals.
For every integer $p\ge1$, we have
$p+m/2\ge2$ and $p+b\ge p+1/2 \ge3/2$, so both bounds in
\eqref{eq:fderthm3} can be integrated using
$\int_0^\infty e^{-A/t}t^{-a}\,dt=\Gamma(a-1)A^{1-a}$,
$A>0$, $a>1$. Multiplying by $|\hat{y}|^{q-s+\gamma}$ and
absorbing its power into the Gaussian in the first bound together with \eqref{eq:etarelation} gives
\begin{align} \label{eq:fintthm3}
\nonumber
        |\hat{y}|^{q-s+\gamma}\int_0^n |f^{(p)}(t)|\,dt
&\le
    \frac{C}{n^{(\eta-m+2)/2}\,\hat{x}^{2p+m-2}}
+
    \frac{C}{n^{m/2}\,|\hat{y}|^{\eta-m+2p}}\\
&\le
    \frac{C}{\sqrt{n}\,M^{\eta+2p-1}},
\end{align}
where we used $\eta-m+1=2j-q+h\ge0$ in the last inequality.

Since $f$ and all its derivatives vanish at $t=0$ and $t=n$,
the Euler--Maclaurin formula
(see, for example, Gel'fond's book \cite{Gelfond}) and
\eqref{eq:fintthm3} with $p=r+1$ give
\begin{align} \label{eq:Snab_int}
        \hat{y}^{q-s+\gamma}\sum_{k=1}^{n-1} f(k)
=
        \hat{y}^{q-s+\gamma}\int_0^n f(t)\,dt
+
    O\left(\frac{1}{\sqrt{n}\,M^{r+1}}\right).
\end{align}

Substituting $t=nw$ and using \eqref{eq:etarelation}, we get
\begin{align*}
    \hat{y}^{q-s+\gamma}\int_0^n f(t)\,dt
=
    \frac{e^{-\frac{(\hat{x}+|\hat{y}|)^2}{2n}}}
    {n^{\eta/2}}
    \left(\frac{\hat{y}}{\sqrt{n}}\right)^{q-s+\gamma}
    J^{(h)}_{m,b}\!\left(\frac{\hat{x}}{\sqrt{n}},
                         \frac{|\hat{y}|}{\sqrt{n}}\right),
\end{align*}
where, for $u, v > 0$,
\begin{align} \label{eq:Jdef}
    J^{(h)}_{m,b}(u,v)
:=
    e^{\frac{(u+v)^2}{2}}
    \int_0^1
        \frac{F^{(h)}_m(u/\sqrt{w})\,
              e^{-\frac{u^2}{2w} - \frac{v^2}{2(1-w)}}}
             {w^{m/2}(1-w)^{b}}\,
    dw.
\end{align}
Recalling that $\eta=2j-q+m+h-1$, we obtain, for $x\ge1, y\neq 0$,
\begin{multline} \label{eq:Tjqs5:thm3}
    T_{j,q,s}(x,y)
=
    e^{-\frac{(\hat{x}+|\hat{y}|)^2}{2n}}
        \sum_{\gamma + 2\mu \le r-s}
        \sum_{m=2}^{r-h+2}
            \frac{(-1)^{\gamma+\mu}}{\gamma!\, \mu!\, 2^\mu}\,
            \frac{(\hat{y}/\sqrt{n})^{q-s+\gamma}}{n^{\eta/2}}
            J^{(h)}_{m,b}\!\left(\frac{\hat{x}}{\sqrt{n}},
                                 \frac{|\hat{y}|}{\sqrt{n}}\right)\\
+
    o\left(\frac{n}{M^{r+3}}\right).
\end{multline}

Substituting \eqref{eq:Tjqs5:thm3} into \eqref{eq:pxydecomp1}
and collecting the remaining terms with the same $\eta$, and omitting terms with $\eta \ge r+2$,
we obtain, for $x\ge1, y \neq 0$,
\begin{align} \label{eq:Pnexpansionthm3}
    \mathcal{P}_n(x,y)
=
    e^{-\frac{(\hat{x}+|\hat{y}|)^2}{2n}}
    \sum_{\eta=1}^{r+1}
        \frac{B_\eta\!\left(\frac{\hat{x}}{\sqrt{n}},
                             \frac{\hat{y}}{\sqrt{n}}\right)}
             {n^{\eta/2}}
+
    o\left(\frac{n}{M^{r+3}}\right),
\end{align}
where, for $u>0$ and $v\ne0$,
\begin{align} \label{eq:Betathm3}
    B_\eta(u,v)
=
    \sum_{2j-q+m+h-1 = \eta}
        a_{q,j}\binom{q}{s}
            \frac{(-1)^{\gamma+\mu}}{\gamma!\, \mu!\, 2^\mu}\,
            v^{q-s+\gamma}\,
            J^{(h)}_{m,b}(u,|v|),
\end{align}
the sum being over nonnegative integers $j, q, s, \gamma, \mu, m$
with
\begin{align*}
    s \le q \le \lf 3j/2 \rf,
    \qquad m \ge 2,
\end{align*}
and
\begin{align*}
    h = s+\gamma+2\mu,
    \qquad b=j+\gamma+\mu+\frac12.
\end{align*}

By \cite[Theorem VII.3.16]{Petrov}, truncating after the term of
order $n^{-(r+1)/2}$, we have, uniformly in $x,y$,
\begin{align} \label{eq:Petrovthm3}
    \pr(S_n = y-x)
=
    e^{-\frac{(\hat{y}-\hat{x})^2}{2n}}
    \sum_{\nu=1}^{r+1}
        \frac{\widehat{q}_\nu\!\left(
                  \frac{\hat{y}-\hat{x}}{\sqrt{n}}\right)}
             {n^{\nu/2}}
+
    O\left(\frac{1}{n^{(r+2)/2}}\right),
\end{align}
and $n^{-(r+2)/2} \le \sqrt{n}M^{-(r+3)}$. Combining this with
\eqref{eq:Pnexpansionthm3} and \eqref{eq:approxPnxy}, we
arrive at, for $x,y\ge1$,
\begin{align*}
    \pr(x+S_n = y, \tau_x > n)
=
    e^{-\frac{(\hat{y}-\hat{x})^2}{2n}}
    \sum_{\nu=1}^{r+1}
        \frac{\widehat{q}_\nu\!\left(
                  \frac{\hat{y}-\hat{x}}{\sqrt{n}}\right)}
             {n^{\nu/2}}
&-
    e^{-\frac{(\hat{x}+\hat{y})^2}{2n}}
    \sum_{\nu=1}^{r+1}
        \frac{B_\nu\!\left(\frac{\hat{x}}{\sqrt{n}},
                            \frac{\hat{y}}{\sqrt{n}}\right)}
             {n^{\nu/2}}\\
&\hspace{3cm}+
    o\left(\frac{n}{M^{r+3}}\right).
\end{align*}
This proves the expansion for $x,y\ge1$. It remains to show that
the restrictions of $B_\nu$ to $u,v>0$ are polynomials and that
$B_1=\widehat{q}_1$ there.

We first use the integral representation to identify the possible
negative powers of $u$ and $v$. Comparing with the unconditional
expansion for $y<0$, we will show that the coefficients of negative
powers of $u$ vanish. Time reversal then excludes negative
powers of $v$.

Throughout the rest of the proof, $u,v>0$.
Consider a standard Brownian motion started at $u$.
We decompose its density at $-v$ at time $1$
according to the first hitting time of zero. Using its density
\cite[Theorem~2.35]{MP10} and the strong Markov property, we obtain
\begin{align} \label{eq:baseintegralthm3}
    \int_0^1
        \frac{e^{-u^2/(2w)-v^2/(2(1-w))}}
             {w^{3/2}(1-w)^{1/2}}\,dw
=
    \frac{\sqrt{2\pi}}u e^{-(u+v)^2/2}
    =:H(u,v).
\end{align}
Define
\begin{align*}
    D_u:=-2\partial_{u^2}=-\frac1u\partial_u,
    \qquad D_v:=-2\partial_{v^2}=-\frac1v\partial_v.
\end{align*}
We differentiate with respect to $u^2$ and $v^2$ so that
$D_u$ and $D_v$ multiply the exponential under the integral
by $1/w$ and $1/(1-w)$, respectively.
Differentiating \eqref{eq:baseintegralthm3} therefore gives
\begin{align*}
    \int_0^1
        \frac{e^{-u^2/(2w)-v^2/(2(1-w))}}
             {w^{3/2+p}(1-w)^{1/2+k}}\,dw
    =D_u^pD_v^k H(u,v),
    \qquad p,k\in\mathbb N_0.
\end{align*}
Each derivative on the right is $e^{-(u+v)^2/2}$ times a
polynomial in $u,v,1/u,1/v$.

By \eqref{eq:Jdef} the monomial of degree $i$ in $F_m^{(h)}$ contributes to
$e^{-(u+v)^2/2}J^{(h)}_{m,b}(u,v)$ a constant multiple of
\begin{align*}
    u^iD_u^pD_v^kH(u,v),\qquad
    p=\frac{m+i-3}{2},\qquad
    k=b-\frac12=j+\gamma+\mu.
\end{align*}
These are nonnegative integers, since $m+i$ is odd and at least
$3$ by Lemma~\ref{lem:momentsoflocprob2}. The Gaussian factor
cancels the prefactor in \eqref{eq:Jdef}, so $J^{(h)}_{m,b}(u,v)$
is a polynomial in $u,v,1/u,1/v$. By \eqref{eq:Betathm3}, the
same holds for both $B_\eta(u,v)$ and $B_\eta(u,-v)$.

For $y<0$, the walk must have exited before or at time $n$,
so we can compare with the unconditional expansion.
We apply \eqref{eq:Pnexpansionthm3} in this range and put
$u=\hat{x}/\sqrt{n}$, $v=-\hat{y}/\sqrt{n}$.
For $u,v\in(1,2)$, the remainder $o(n/M^{r+3})$ is
$o(n^{-(r+1)/2})$. Clearly,
\begin{align*}
    \pr(x+S_n=y,\tau_x=n)
\le
    \pr(X_1<y)
=
    O\left(n^{-(r+3)/2}\right).
\end{align*}
Consequently, \eqref{eq:approxPnxy} gives
$\mathcal{P}_n(x,y)=\pr(S_n=y-x)+o(n^{-(r+1)/2})$.
The expansions \eqref{eq:Pnexpansionthm3} and \eqref{eq:Petrovthm3}
have the same Gaussian factor
$e^{-(u+v)^2/2}$, which is bounded away from zero in this range.
Cancelling it, we obtain, uniformly in this range,
\begin{align} \label{eq:Bminuscomparisonthm3}
    \sum_{\nu=1}^{r+1}\frac{1}{n^{\nu/2}}
    \left[
        B_\nu(u,-v)-\widehat{q}_\nu(-v-u)
    \right]
=o\left(n^{-(r+1)/2}\right).
\end{align}

Fix $u,v\in(1,2)$ and take $x=\lf u\sigma\sqrt{n}\rf$,
$y=-\lf v\sigma\sqrt{n}\rf$. Multiplying
\eqref{eq:Bminuscomparisonthm3} by $\sqrt{n}$ and passing to the
limit gives
$B_1(u,-v)=\widehat{q}_1(-v-u)$ by continuity.
After clearing denominators, this is a polynomial identity on
$(1,2)^2$, hence it holds for all $u,v>0$. We can therefore
subtract the first coefficients before comparing the next terms.
Repeating this argument gives
\begin{align} \label{eq:Bminusidentity}
    B_\nu(u,-v) = \widehat{q}_\nu(-v-u),
    \qquad u,v>0,\ 1 \le \nu \le r+1.
\end{align}

To transfer this cancellation to $B_\eta(u,v)$, we collect the
derivatives obtained from \eqref{eq:Betathm3} into the finite
operator
\begin{align*}
    \mathcal D_\eta
    =\sum_{i,d,p,k}c^{(\eta)}_{i,d,p,k}\,
        u^iv^d D_u^pD_v^k.
\end{align*}
Here $i,p,k$ are as above, $d=q-s+\gamma\ge0$, and the
coefficients $c^{(\eta)}_{i,d,p,k}$ are constants. Then
\begin{align} \label{eq:Doperatorthm3}
    e^{-(u+v)^2/2}B_\eta(u,v)=\mathcal D_\eta H(u,v).
\end{align}
Replacing the external factors $v^d$ by $(-v)^d$ gives
$B_\eta(u,-v)$ instead. Hence \eqref{eq:Bminusidentity} yields
\begin{align} \label{eq:Bnegativeoperatorthm3}
    \sum_{i,d,p,k}c^{(\eta)}_{i,d,p,k}\,
        u^i(-v)^d D_u^pD_v^k H(u,v)
    =e^{-(u+v)^2/2}\widehat q_\eta(-v-u).
\end{align}
We use the explicit expression for $H$ also at negative values
of its second argument. After cancelling the Gaussian factor
and clearing denominators, \eqref{eq:Bnegativeoperatorthm3}
is a polynomial identity, so we may replace $v$ by $-v$.
Since $D_v[f(-v)]=(D_vf)(-v)$, restoring the Gaussian factor
gives
\begin{align} \label{eq:Hreflectedthm3}
    \mathcal D_\eta[H(u,-v)]
    =e^{-(u-v)^2/2}\widehat q_\eta(v-u),
    \qquad u,v>0.
\end{align}

We now subtract \eqref{eq:Hreflectedthm3} from
\eqref{eq:Doperatorthm3} to show that $B_\eta(u,v)$ remains
bounded as $u\to0$. We have
\begin{align} \label{eq:Bregularthm3}
    e^{-(u+v)^2/2}B_\eta(u,v)
    &=e^{-(u-v)^2/2}\widehat q_\eta(v-u)
    \nonumber\\
    &\quad+\mathcal D_\eta[H(u,v)-H(u,-v)].
\end{align}
The difference
\begin{align*}
    H(u,v)-H(u,-v)
    =-2\sqrt{2\pi}\,e^{-(u^2+v^2)/2}\frac{\sinh(uv)}u,
\end{align*}
has a convergent power series in $u^2$ at zero, since
$\sinh(uv)/u=v+u^2v^3/3!+\cdots$.
The operators $D_u=-2\partial_{u^2}$ and $D_v$ preserve this
property for $v>0$. Multiplication by $u^iv^d$, with $i,d\ge0$,
then leaves a function bounded at $u=0$ for each fixed $v>0$.
Thus the right-hand side of \eqref{eq:Bregularthm3}
is bounded as $u\to0$. The Gaussian factor on the left tends
to $e^{-v^2/2}>0$, so $B_\eta(u,v)$ is bounded as well.
Its polynomial representation in $u,v,1/u,1/v$ therefore
contains no negative powers of $u$.

The negative powers of $v$ are excluded by time reversal:
$$\pr(x+S_n = y, \tau_x > n) = \pr(y+S^*_n = x, \tau^*_y > n).$$
The reflected walk $S^*_n$ satisfies the same moment assumption.
The expansion proved above applies to this walk with
coefficients $B^*_\nu$. The unconditional terms agree because
$\pr(S^*_n=x-y)=\pr(S_n=y-x)$. Comparing the coefficients as in
the proof of \eqref{eq:Bminusidentity}, we obtain
\begin{align} \label{eq:Bsymmetrythm3}
    B_\nu(u,v)=B^*_\nu(v,u),
    \qquad u,v>0,\ 1\le\nu\le r+1.
\end{align}
The argument above applies also to the reflected walk, so
$B^*_\nu$ has no negative powers of its first argument.
By \eqref{eq:Bsymmetrythm3}, $B_\nu$ has no negative powers of
$v$ either.

For $\eta=1$, the indices in \eqref{eq:Betathm3} satisfy
$q=s=\gamma=0$, so $B_1(u,v)=B_1(u,-v)$. Since
$\widehat{q}_1$ is constant, \eqref{eq:Bminusidentity} gives
$B_1=\widehat{q}_1$ also for $u,v>0$.

This completes the proof of Theorem~\ref{thm:3}.

\section*{Statement on the use of generative AI}

All main ideas of the proof were created without using generative AI. 
It was widely used for proofreading the text and to simplify the technical computations. 
The proof of Theorem~\ref{thm:3} follows the same scheme as in the proof of 
Theorem~\ref{thm:2} and AI was used to work out details such as handling small 
times $k$ with $k \le C x^2/ \log n$ and proving that the functions appearing in 
the decomposition do not contain negative powers of variables and so they are polynomial. 

All suggestions produced by the AI were critically examined and
independently verified by the author.

\bibliographystyle{abbrv}
\bibliography{references}
\end{document}